\documentclass [preprint,12pt]{elsarticle}
\usepackage{amssymb}
\journal{Journal of Mathematical Chemistry}

\begin{document}


\title
{Isochronicity and limit cycle oscillation in chemical systems}

\author{Sandip Saha }
\address{S N Bose National Centre For Basic Sciences\\Block-JD, 
Sector-III, Salt Lake, Kolkata-700106, India.}

\author{Gautam Gangopadhyay}

\address{S N Bose National Centre For Basic Sciences\\Block-JD, 
Sector-III, Salt Lake, Kolkata-700106, India.}

\ead{gautam@bose.res.in}

\date{\today}

\begin{abstract}

Chemical oscillation is an interesting nonlinear dynamical phenomenon which arises due to complex stability condition of the 
steady state of a reaction far away from equilibrium which is usually characterised by a periodic attractor or a limit cycle around an 
interior stationary point. In this context $Lienard$ equation is specifically used in the study of nonlinear dynamical properties of an open 
system which can be utilized to obtain the condition of limit cycle.
 In conjunction with the property of  limit cycle oscillation, here we have shown  the condition for isochronicity for different chemical oscillators with the help of renormalisation group method with 
multiple time scale analysis  from a $Lienard$ system. When two variable open system of equations are transformed into a $Lienard$
system of equation the condition for limit cycle and isochronicity can be stated in a
unified way. 
 For any such nonlinear oscillator we have shown the route  of a dynamical transformation of a limit cycle oscillation to a periodic orbit of centre type depending on the parameters of the system. 

\end{abstract}


\maketitle
\section{Introduction}

To bypass the exact solution of  the non-linear dynamical systems the general trend is  to resort to a geometrical approach coupled with 
tools of analysis\cite{nayfeh,birkhoff,andronov,arnold,smale,ruelle}. For example, in many body dynamics thus people are more interested 
in finding wheather the system is going  to be stable forever or will one or few  bodies dissociate from the rest. Most non-linear dynamical 
systems are handled by various
perturbative approaches or asymptotic analyses where the perturbation theory usually confers to a collection of iterative methods 
for the systematic analysis of global behaviour of differential equations. As ordinary perturbation theory often  fails due 
to non-convergence\cite{nayfeh,birkhoff,bender} of the series so in order to extract information from perturbation theory
 there is a need to develop proper techniques to tackle the summation of otherwise divergent series. In this circumstance one has to look for  various singular perturbation techniques. All these  methods basically demands that at every order of perturbation the so called secular or divergent terms arising out of straightforward application of perturbation theory be removed. It has proved to be a successful tool in finding approximate solutions to weakly non-linear differential equations with finite oscillatory period. In the perturbative renormalization treatment of dynamical systems the  amplitude and phase  of the oscillation get renormalized\cite{wilson,holmesbook}. When the approximate solution is expressed in terms of these renormalized amplitude and phase then the perturbative series is uniformly valid and does not have any secular term\cite{chen1,chen2}. In the traditional  renormalization Group($RG$) approach \cite{wilson} in field  theory the order parameters designated by  various coupling constants play similar role as amplitude and phase in the case of oscillatory dynamical system. $RG$ method also deals with the problem of isochronous centres characteristic of a family of initial condition dependent periodic orbits\cite{len0,len1,len2,len3,len3.5} surrounding a critical point. Isochronicity is a widely studied subject not only for its relation with stability theory and bifurcation theory\cite{bender,murraynld,murraymbio,sabatini} but also in the study of bifurcations of critical points leading to limit cycles\cite{len4} and isochronous systems\cite{len1,len2,len3}.

Various kinds of periodic trajectories in phase space can be found and most striking example is a self-sustained oscillation or limit cycle in the system. A centre refers to a family of initial-condition-dependent periodic orbits surrounding a point. While centres can exist in both linear and non-linear systems, limit cycles can occur only in non-linear systems. The most important kind of limit cycle is the stable limit cycle where all nearby trajectories spiral towards the isolated orbit. Existence of a stable limit cycle in a dynamical system means there exists self-sustained oscillation in the system. Dynamical systems capable of having limit cycle oscillations are very important from the point of view of modelling real-world systems which exhibit self-sustained oscillations. Some examples of such phenomena from nature include heart beating\cite{bab,pierson}, oscillations in body temperature\cite{refi,jewett}, random eye movement oscillations during sleep\cite{mac}, hormone secretion\cite{glass,gold}, chemical reactions that oscillate spontaneously\cite{prigogine,richard,schnaken,kuramoto} etc. There can be other varieties of limit cycles as well: (i) the unstable limit cycle, where all neighbouring trajectories tend to move away from the isolated orbit in phase space and (ii) semi-stable limit cycle, where trajectories are attracted to the limit cycle on one side and repelled on the other. In a sense centres can be thought of as neutrally stable limit cycles. Although a lot of work has been performed in finding ways to determine if a system has a limit cycle, surprisingly a little is known about how to find this and still it remains a highly active area of research\cite{goldbook,epstein, strogatz,sarkar}.

Chemical oscillation\cite{grayscott,epstein} is an interesting non-linear dynamical phenomenon which arises due to intrinsic instability of the non-equilibrium steady state of a reaction under far away from equilibrium condition\cite{prigogine}. Experimentally such open systems like, Bray\cite{bray}, BZ\cite{bz1,bz2,bz3,lav1,lav2}  and glycolytic reactions\cite{gold,gly1,glyscott,gly2,gly3} are studied extensively in a continuously
flowing stirred tank reactor and the nature of the oscillatory kinetics of two 
intermediates gave  reliable dynamical models of limit cycle. The self-sustained chemical oscillations\cite{prigogine,grayscott,epstein} are also regularly observed in biological world to maintain a cyclic steady state e.g., glycolytic oscillations\cite{glyscott,goldbook,gly2,gly3}, calcium  oscillations\cite{cal}, cell division\cite{cel}, circadian oscillation\cite{cir} and others\cite{epstein}. The generic features of such diverse nature of non-linear chemical oscillations are due to auto-catalysis and various feedback mechanisms into the system which are basically controlled by a few slowest time scales of the overall process. The coupled dynamics of the system can be described by two intermediate concentrations or population variables  characterized by the occurrence of a limit cycle when the motion is visualized on a phase plane.
In particular, a procedure of the reduction of chemical cubic equations to the form of a second order differential equation with co-efficients which allows for the limit cycle analysis so called Rayleigh oscillator has been proposed for a first time in the article by Lavrova et al\cite{lav1,lav2}. Its further development to a more general case, the Lienard oscillator was given by Ghosh and Ray\cite{len4}.  Here we consider that a class of arbitrary, autonomous kinetic equations in two variables describing chemical oscillations can be cast into the form of a $Lienard$ oscillator\cite{len0,len1,len2,len3,len3.5,len4}. It is characterised by  the non-linear  forcing and damping coefficients which can control the limit cycle behavior.

 Although nonlinear oscillators got a lot of attention over the years 
in the field of dynamical systems but there is no  straightforward way of distinguishing between
limit cycle and isochronous orbit of the system with their very different kinds of solutions. 
Here our effort is  to study the isochronous systems by analyzing the
behaviour around a centre by  $RG$ method\cite{chen1,chen2} and to find the so called
periodic functions i.e. conditions under which a system becomes isochronous\cite{len3,len3.5}. 
 Our method of distinguishing center and limit cycle behaviours are numerically analyzed here in terms of the parameters of the chemical oscillators.
  When the two dimensional kinetic equations are transformed into a $Lienard$
system of equation we would like to find here the relation between the limit cycle and isochronicity
 in an open dynamical system. As the chemical oscillators are standard real experimental models the theory is verified here in various systems.
 
 In section 2, we have briefly reviewed the method of reduction of kinetic equation into $Lienard$ form to find the condition for limit cycle. Isochronicity for  $Lienard$ System is described in section 3. In section 4, we have shown the examples, (4.1) modified Brusselator model, (4.2) Glycolytic oscillator and (4.3) $van$ $der$ $Pol$ type oscillator to analyze the behaviour of limit cycle and isochronicity. The paper is concluded in section 5.

\section{Reduction of Kinetic Equation into $Lienard$ Form: Conditions for limit cycle}

Let us consider a two dimensional set of autonomous kinetic equations for open system.  
Here our purpose is to review the condition for limit cycle by casting the two dynamical equations into a form of $Lienard$ oscillator\cite{strogatz,len0,len1,len2,len3,len4}.
Following the analysis of Ghosh and Ray\cite{len4} we  consider a system of differential equations 

$$\frac{dx}{dt} = a_0+a_1x+a_2y+f(x,y),$$
\begin{equation}
\frac{dy}{dt} = b_0+b_1x+b_2y+g(x,y),
\label{eq1}
\end{equation}
where $x$ and $y$ are populations of two intermediate species of a dynamical process with $a_0,a_1,a_2,b_0,b_1,b_2$ are all real parameters expressed in terms of the kinetic constants with $f(x,y)$ and $g(x,y)$ are  non-linear functions.

Then  writing the equations in terms of a new pair, $(z,u)$ as
  $$u=\alpha_0+\alpha_1x+\alpha_2y ,$$
\begin{equation}
z=\beta_0+\beta_1x+\beta_2y, 
\label{eq2}
\end{equation}
where $\alpha_0,\alpha_1,\alpha_2$ and $\beta_0,\beta_1,\beta_2$
are constants  expressed in terms of $a_i$ and $b_i$.
Choosing $u$ and $z$ in such a way that, 
\begin{equation}
\frac{dz}{dt}=u,
\label{eq4}
\end{equation}
and differentiating (\ref{eq4})  w.r.t. t, we get

$$\ddot{z}=\dot{u}=\alpha_1\dot{x}+\alpha_2\dot{y}=
\alpha_1 \lbrace a_0+a_1 L(z,\dot{z})+a_2 M(z,\dot{z})+\varphi(z,\dot{z})\rbrace$$
\begin{equation}
\hspace{3.3 cm}+\alpha_2 \lbrace b_0+b_1 L(z,\dot{z})+b_2 M(z,\dot{z})+\phi(z,\dot{z})\rbrace
\label{eq5}
\end{equation}
where 
$$L(z,\dot{z})=c_1 z+c_2\dot{z}+c_L$$ 
with  $$c_1=-\frac{\alpha_2}{\alpha_1 \beta_2-\alpha_2 \beta_1},c_2=\frac{\beta_2}{\alpha_1 \beta_2-\alpha_2 \beta_1},c_L=\frac{\alpha_2 \beta_0-\alpha_0 \beta_2}{\alpha_1 \beta_2-\alpha_2 \beta_1}.$$ and
$$M(z,\dot{z})=c_3 z+c_4\dot{z}+c_M$$ with $$c_3=-\frac{\alpha_1}{\alpha_2 \beta_1-\alpha_1 \beta_2},c_4=\frac{\beta_1}{\alpha_2 \beta_1-\alpha_1 \beta_2},c_M=\frac{\alpha_0 \beta_1-\alpha_1 \beta_0}{\alpha_2 \beta_1-\alpha_1 \beta_2}.$$

Next we consider $c_L$ and $c_M$ to be negligibly small. It is trivial to assume that  both the numerators will not exactly vanish as  $\alpha_2 \beta_0=\alpha_0 \beta_2$ and $\alpha_0 \beta_1=\alpha_1 \beta_0$ together giving  $\alpha_2 \beta_1=\alpha_1 \beta_2$ which makes all constants $c_i$ and the system to be undefined. Here it is performed  by choosing the  ratio of numerator and denominator for the constants $c_L$ and $c_M$ are very small. Subsequently  we define $L(z,\dot{z})$ and $M(z,\dot{z})$ by ignoring the small values of $c_L$ and $c_M$, respectively. 

Now taking the functions $\varphi$ and $\phi$  as power series i.e. 
\begin{equation}
\varphi(z,\dot{z})=\sum_{n,m=0}^{\infty} \varphi_{nm} z^n \dot{z}^m \hspace{1cm}\phi(z,\dot{z})=\sum_{n,m=0}^{\infty} \phi_{nm} z^n \dot{z}^m
\label{eq7}
\end{equation} 
in equation (\ref{eq5})  and using  (\ref{eq7}) one finds


\begin{equation}
\ddot{z}=A_{00}+A_{10} z+A_{01} \dot{z}+\sum_{n>1}^{\infty} A_{n0} z^n+\sum_{m>1}^{\infty} A_{0m} \dot{z}^m+\sum_{n,m\geqslant1}^{\infty} A_{nm} z^n \dot{z}^m ,
\label{eq9}
\end{equation}
where, $A_{00}= \alpha_1 a_0+\alpha_2 b_0+\alpha_1 \varphi_{00}+\alpha_2 \phi_{00}$, $A_{10}=\alpha_1(a_1 c_1+a_2 c_3)+\alpha_2 (b_1 c_1+b_2 c_3)+\alpha_1 \varphi_{10}+\alpha_2 \phi_{10}$, $A_{01}=\alpha_1(a_1 c_2+a_2 c_4)+\alpha_2 (b_1 c_2+b_2 c_4)+\alpha_1 \varphi_{01}+\alpha_2 \phi_{01}$, $A_{n0}=\alpha_1 \varphi_{n0}+\alpha_2 \phi_{n0}$, $A_{0m}=\alpha_1 \varphi_{0m}+\alpha_2 \phi_{0m}$ and $A_{nm}=\alpha_1 \varphi_{nm}+\alpha_2 \phi_{nm},\forall n,m\geqslant 1. $

Now, for the steady state, $z=z_s$, both $\dot{z}$ and $\ddot{z}$ vanish and the fixed points follow the condition
\begin{equation}
A_{00}+A_{10} z_s+\sum_{n>1}^{\infty} A_{n0} z_s^n = 0.
\label{eq10}
\end{equation}

The equation for deviation from the  stationary point from $z$ i.e. $\xi(=z-z_s)$ follows from equation (\ref{eq9}) as,

\begin{equation}
\ddot{\xi}+F(\xi,\dot{\xi}) \dot{\xi} +G(\xi)=0,
\label{eq12}
\end{equation}
where the functions $F(\xi,\dot{\xi})$ and $G(\xi)$ are given by
$$
F(\xi,\dot{\xi})=-[A_{01}+\sum_{m>1}^{\infty} A_{0m} \dot{\xi}^{m-1}+\sum_{n,m\geq1}^{\infty} A_{nm} (\xi+z_s)^n \dot{\xi}^{m-1}]$$
\begin{equation}
G(\xi)=-[A_{00}+A_{10} (\xi+z_s)+\sum_{n>1}^{\infty} A_{n0} (\xi+z_s)^n ].
\end{equation}

Equation (\ref{eq12}) is a well known form of generalised $Lienard$ equation if the damping force, $F(\xi,\dot{\xi})$ and the restoring force, $G(\xi)$  satisfy the usual regularity conditions as given in Strogatz\cite{strogatz} page-210.





So, the condition for existence of having a stable limit cycle of the above described $Lienard$ system should satisfy $F(0,0)<0$ i.e.
\begin{equation}
-[A_{01}+ \sum_{n\geqslant1}^{\infty} A_{n1} z_s^n ]<0.
\end{equation}

\section{Isochronicity for  $Lienard$ System}

For a given two dimensional non-linear dynamical system of equations, in general, can be cast into  $Lienard$ system.
From the $Lienard$ system one can set up a perturbation theory around the closed orbit of the centre\cite{len1,len2,len3,len3.5}. The orbit is characterised by two constants, the amplitude, $A$ and the phase, $\theta$, fixed by the two initial conditions. However, the perturbation theory  most likely diverges due to the presence of secular terms as the separation of time scale becomes large. Two renormalisation constants have to be introduced to absorb these 
divergence. 
 The renormalisation constants appear in terms of an arbitrary time, say $\tau$, which serves to fix the new initial condition which makes the 
amplitude and the phase then dependent on $\tau$.
 The value of $x$ at $t$ cannot depend on  where one sets the initial condition and hence $(\frac{\partial x}{\partial \tau})_t =0$, which is 
the flow equation. This must give $\frac{d A}{d \tau} =p(A)$ and $\frac{d \theta}{d \tau} =q(A)$. If the system is of centre type then the initial condition sets the amplitude of motion and hence $\frac{d A}{d \tau} =0$. The phase flow equation on the other hand normally furnishes the non-linear correction to the frequency. However, for an isochronous centre there can be no correction to the frequency and hence $p(A)=0$ and $q(A)=0$ identically.

Use of above technique can serve in  differentiating between oscillatory dynamics of a  centre type or limit cycle. The centre type oscillation consists of a continuous family of closed orbits in phase space, each orbit being determined by its own initial condition. This implies the amplitude, $A$ is fixed once the initial condition is set\cite{len3}.
Our objective is to derive the condition for isochronicity\cite{len3,len3.5} from  a $Lienard$ oscillatory equation and finally the relation between the condition for being a stable limit cycle and the mutual relation.

First we find  a simple $Lienard$ oscillatory system from (\ref{eq12}) by taking some special order of the power series in which $n,m$  contribute starting from $0$ to atmost $2$ and from that we may get a polynomial of highest degree atmost 3 in the damping force function for the $Lienard$ system.
Using the above assumption, one can obtain, 

$$F(\xi,\dot{\xi})=-[A_{01}+\sum_{m>1}^{2} A_{0m} \dot{\xi}^{m-1}+\sum_{n,m\geq1}^{2} A_{nm} (\xi+z_s)^n \dot{\xi}^{m-1}]$$
$\Rightarrow$
$$F(\xi,\dot{\xi})=-[A_{01}+A_{02} \dot{\xi}+A_{11} {\xi}+A_{11} z_s+A_{12} \xi \dot{\xi}+A_{12} z_s \dot{\xi}+A_{21} \xi^2$$
\begin{equation}
\hspace{1.5 cm}+2A_{21} \xi z_s+A_{21} z_s^2+A_{22} \xi^2 \dot{\xi}+2A_{22} \xi z_s \dot{\xi}+A_{22} z_s^2 \dot{\xi}]
\end{equation}
and
$$G(\xi)=-[A_{00}+A_{10} \xi+A_{10} z_s+A_{20} \xi^2+2A_{20} \xi z_s+A_{20} z_s^2].$$

The condition for being stable limit cycle, $F(0,0)<0$ gives,
\begin{equation}
\hspace{3 cm}-[A_{01}+A_{11} z_s+A_{21} z_s^2]<0
\label{eq16}
\end{equation}
and (\ref{eq10}) gives
\begin{equation}
\hspace{3.5 cm}A_{00}+A_{10} z_s+A_{20} z_s^2 = 0
\end{equation}
and therefore, $$G(\xi)=-[A_{10} \xi+A_{20} \xi^2+2A_{20} \xi z_s].$$
As $G(\xi)$ must satisfy $G(\xi)=0$ at $\xi=0$, 
thus the $Lienard$ system becomes, 
$$\ddot{\xi}-[A_{01}+A_{02} \dot{\xi}+A_{11} {\xi}+A_{11} z_s+A_{12} \xi \dot{\xi}+A_{12} z_s \dot{\xi}+A_{21} \xi^2+2A_{21} \xi z_s+A_{21} z_s^2$$
\begin{equation}
+A_{22} \xi^2 \dot{\xi}+2A_{22} \xi z_s \dot{\xi}+A_{22} z_s^2 \dot{\xi}]\dot{\xi}
-[A_{10} \xi+A_{20} \xi^2+2A_{20} \xi z_s]=0.
\end{equation}

Now if we set the condition other than stable limit cycle i.e. $F(0,0)=0$, for example, taking (\ref{eq16}) as zero for some values of the parameters, then above equation on simplifying becomes, 

$$\ddot{\xi}+\omega^2 \xi=A_{02} \dot{\xi}^2+A_{11} {\xi} \dot{\xi}+A_{12} \xi \dot{\xi}^2+A_{12} z_s \dot{\xi}^2+A_{21} \xi^2 \dot{\xi}+2A_{21} \xi z_s \dot{\xi}$$
\begin{equation}
\hspace{1.5 cm}+A_{22} \xi^2 \dot{\xi}^2+2A_{22} \xi z_s \dot{\xi}^2+A_{22} z_s^2 \dot{\xi}^2+A_{20} \xi^2
\label{eq19}
\end{equation}
where $\omega^2=-2A_{20} z_s-A_{10}$ must be $+ve$.  Since $\omega$ is a real quantity and for  $\omega^2<0$ then $G(\xi)$  violate its property.

For book keeping purpose we introduce a positive $\lambda$, with $\lambda<<1$ and using on (\ref{eq19})
and  discarding higher orders of $\lambda$ we get,

\begin{equation}
\ddot{\xi}+\omega^2 \xi=\lambda A_{20} \xi^2+\lambda [A_{02}+A_{22} z_s^2+A_{12} z_s] \dot{\xi}^2+\lambda [A_{11}+2A_{21} z_s] {\xi} \dot{\xi}+O(\lambda^2).
\label{eq20}
\end{equation}

After that using Renormalisation  Group ($RG$) technique\cite{len3}, let  us take a perturbation solution of $\xi=\xi_0+\lambda \xi_1+\lambda^2 \xi_2+\lambda^3 \xi_3+\cdots\cdots\cdots$ , i.e. $\xi=\xi_0+\lambda \xi_1+O(\lambda^2)$, to get an approximate solution of (\ref{eq20}). 
So, putting $\xi$ and after simplifying (on neglecting $O(\lambda^2)$), we get,
$$(\ddot{\xi_0}+\ddot{\lambda \xi_1})+\omega^2 (\xi_0+\lambda \xi_1)=\lambda A_{20} \xi_0^2+\lambda(A_{02}+A_{12} z_s+A_{22} z_s^2) \dot{\xi_0}^2$$
\begin{equation}
\hspace{4.5cm}+\lambda(A_{11}+2 A_{21} z_s) \xi_0 \dot{\xi_0}+O(\lambda^2).
\end{equation}

Comparing the co-efficient of $\lambda^0$, $\lambda^1$ of both sides we get,
\begin{equation}
\lambda^0 : \ddot{\xi_0}+\omega^2 \xi_0=0
\end{equation}
\begin{equation}
\lambda^1 : \ddot{\xi_1}+\omega^2 \xi_1=A_{20} \xi_0^2+(A_{02}+A_{12} z_s+A_{22} z_s^2) \dot{\xi_0}^2+(A_{11}+2 A_{21} z_s) \xi_0 \dot{\xi_0}.
\end{equation}

If we take higher order terms then it must be included within $O(\lambda^2)$ and we simply neglect here $O(\lambda^2)$. Let us set an initial condition $\xi(t)=A$ and $\dot{\xi(t)}=0$ at $t=t_0$ with $t_0$ being the initial time, then by 
comparing $\lambda$ as  previously we  get $\xi_0=A$ and $\xi_i=0, \forall i>0$  
along with  $\dot{\xi_i}=0, \forall i\geq0$ at $t=t_0$.

Thus after solving above equations we get $\xi_0(t)$ and $\xi_1(t)$ as,

$$\xi_0(t)=A \cos \omega (t-t_0)$$
$\xi_1(t)=-[\frac{A^2 A_{20}}{3 \omega^2}+\frac{2 A^2 (A_{02}+A_{12} z_s+A_{22} z_s^2)}{3}]\cos\omega(t-t_0) - \frac{A^2 (A_{11}+2 A_{21} z_s) \sin\omega(t-t_0)}{3 \omega}
$

\hspace{1 cm}$+\frac{A^2 A_{20}}{2} [\frac{1}{\omega^2}-\frac{cos 2\omega(t-t_0)}{3 \omega^2}]
+\frac{(A_{02}+A_{12} z_s+A_{22} z_s^2) \omega^2 A^2}{2} [\frac{1}{\omega^2}+\frac{cos 2\omega(t-t_0)}{3 \omega^2}]$
\begin{equation}
+\frac{A^2 (A_{11}+2 A_{21} z_s) \sin 2\omega(t-t_0)}{6 \omega}.
\end{equation}  

So the approximate solution of $\xi(t)$ is, 

$\xi(t)=A \cos \omega (t-t_0)-\lambda [\frac{A^2 A_{20}}{3 \omega^2}+\frac{2 A^2 (A_{02}+A_{12} z_s+A_{22} z_s^2)}{3}]\cos\omega(t-t_0)$

\hspace{1 cm}$-\lambda \frac{A^2 (A_{11}+2 A_{21} z_s) \sin\omega(t-t_0)}{3 \omega}+\lambda \frac{A^2 A_{20}}{2} [\frac{1}{\omega^2}-\frac{cos 2\omega(t-t_0)}{3 \omega^2}]
+\lambda \frac{(A_{02}+A_{12} z_s+A_{22} z_s^2) \omega^2 A^2}{2} $
\begin{equation}
[\frac{1}{\omega^2}+\frac{cos 2\omega(t-t_0)}{3 \omega^2}]+\lambda \frac{A^2 (A_{11}+2 A_{21} z_s) \sin 2\omega(t-t_0)}{6 \omega}
\end{equation}
where $A$ is  the amplitude and $\omega$ is frequency supposing the constant $-\omega t_0 = \theta_0$.

Then $\xi(t)$  becomes,

$\xi(t)=A \cos (\omega t+\theta_0)-\lambda [\frac{A^2 A_{20}}{3 \omega^2}+\frac{2 A^2 (A_{02}+A_{12} z_s+A_{22} z_s^2)}{3}]\cos(\omega t+\theta_0)$

\hspace{1 cm}$-\lambda \frac{A^2 (A_{11}+2 A_{21} z_s) \sin(\omega t+\theta_0)}{3 \omega}+\lambda \frac{A^2 A_{20}}{2} [\frac{1}{\omega^2}-\frac{cos 2(\omega t+\theta_0)}{3 \omega^2}]
+\lambda \frac{(A_{02}+A_{12} z_s+A_{22} z_s^2) \omega^2 A^2}{2} $
\begin{equation}
[\frac{1}{\omega^2}+\frac{cos 2(\omega t+\theta_0)}{3 \omega^2}]+\lambda \frac{A^2 (A_{11}+2 A_{21} z_s) \sin 2(\omega t+\theta_0)}{6 \omega}.
\label{eq26}
\end{equation}

At this point  add another perturbation in the time interval $(t-t_0)$, by splitting $(t-\tau)+(\tau-t_0)$, where $t_0<\tau<t$ and $\tau$ is very close to $t_0$ by defining the interval $(t-\tau)$ as a principal part and the remaining part $(\tau-t_0)$ can be neglected because of smallness. 

Suppose that taking perturbation the time interval, amplitude and phase will be slightly changed from $A$ to $A(\tau)$ and $\theta_0$ to $\theta(\tau)$ . From $RG$ technique the relation between them are $A(\tau)=\frac{A}{Z_1(\tau,t_0)}$ and $\theta(\tau)=\theta_0-Z_2(\tau,t_0)$, where 
\begin{equation}
Z_1(\tau,t_0)=1+\sum_{1}^{\infty} \lambda^n p_n \hspace{1 cm} and \hspace{1 cm} Z_2(\tau,t_0)=0+\sum_{1}^{\infty} \lambda^n q_n.
\end{equation}

Neglecting terms of $O(\lambda^2)$   we get,
\begin{equation}
Z_1(\tau,t_0)=1+\lambda p_1+O(\lambda^2) \hspace{.5 cm} and \hspace{.7 cm} Z_2(\tau,t_0)=\lambda q_1+O(\lambda^2).
\end{equation} 

Now if we put the function $Z_1$ and $Z_2$ as well as $A$ and $\theta_0$ in (\ref{eq26}) and remove the terms which could  led to divergence, we must get either $p_1$ is zero or anything containing $(\tau-t_0)$ and the same for $q_1$ also. But because of the smallness of $(\tau-t_0)$, we can take $p_1$ and $q_1$ approximately to be zero. So, after considering  above, the constants $A$ become $A(\tau)$ and  $\theta_0$ become $\theta(\tau)$, i.e. they become dependent upon the time variable, $\tau$.  Also, if any term   multiplied directly by $(t-t_0)$ in the final solution of $\xi(t)$, then we can convert it into $(t-\tau)$ by neglecting the other part. But here no such terms are directly involved 
in this solution.

So $\xi(t)$ becomes,

$\xi(t)=A(\tau) \cos (\omega t+\theta(\tau))-\lambda [\frac{A^2(\tau) A_{20}}{3 \omega^2}+\frac{2 A^2(\tau) (A_{02}+A_{12} z_s+A_{22} z_s^2)}{3}]\cos(\omega t+\theta(\tau))$

\hspace{1 cm}$-\lambda \frac{A^2(\tau) (A_{11}+2 A_{21} z_s) \sin(\omega t+\theta(\tau))}{3 \omega}+\lambda \frac{A^2(\tau) A_{20}}{2} [\frac{1}{\omega^2}-\frac{cos 2(\omega t+\theta(\tau))}{3 \omega^2}]
+\lambda \frac{(A_{02}+A_{12} z_s+A_{22} z_s^2) \omega^2 A^2(\tau)}{2} $
\begin{equation}
[\frac{1}{\omega^2}+\frac{cos 2(\omega t+\theta(\tau))}{3 \omega^2}]+\lambda \frac{A^2(\tau) (A_{11}+2 A_{21} z_s) \sin 2(\omega t+\theta(\tau))}{6 \omega}.
\end{equation}

Since the final solution cannot depend on the arbitrary time scale, $\tau$, we impose the condition $(\frac{\partial \xi}{\partial \tau})_t=0$ which leads to

\begin{equation}
\frac{dA}{d\tau}=0 \hspace{1cm}\&\hspace{1cm} \frac{d\theta}{d\tau}=0.
\end{equation}
The independence of $\theta$ upon $\tau$ i.e. $\frac{d\theta}{d\tau}=g(A)=0$ gives the condition for isochronicity. If it is non-zero then the $Lienard$ system would not be isochronous. Thus the system will be isochronous for any values of $A$ only when $F(0,0)=0$.

Further, if $\frac{dA}{d\tau}=f(A)$ then we can say there be a limit cycle if $f(A)\neq0$ and the radius of the cycle can be obtained by making $f(A)=0$ if any non-zero $A$  is  found. Otherwise we cannot have any limit cycle because we cannot get any idea about the radii of the cycle. If this type of difficulty comes then we may call this as a centre type. When $\frac{dA}{d\tau}=f(A)=0$ then it is also called centre type.

So, it is seen that, by pushing the limit cycle condition $[F(0,0)<0]$ as zero i.e. $F(0,0)=0$ which is the constant portion  present in the damping force, the $Lienard$ type oscillator transforms into an isochronous oscillator and this is the only condition for being isochronous oscillator. For  $F(0,0)=0$ finally it shows that the $Lienard$ system looses its stability as limit cycle  and becomes a centre type.

\section{Some Chemical oscillator models}

Here we consider a few chemical oscillator models as examples of open system. In open systems there are some inputs and outputs, however, it is possible to attain a steady state depending upon the values of the parameters in addition to dynamical complexities  due to the non-linearities of the system of equations. Inspired by the above analysis of $Lienard$ system we now study some examples to check and verify the above results of limit cycle and isochronicity.

\subsection{Modified Brusselator Model}

The classical $Brusselator$ model\cite{prigogine,epstein,grayscott} is known 
to exhibit kinetics of model trimolecular irreversible reactions which are  based on the vast studies of chemical oscillations\cite{bray,bz1,bz2,bz3} in various systems. The reduction of the Brusselator model in the form of Rayleigh\cite{lav1,lav2} and $Lienard$ form\cite{len4} of differential equations are already published. Here we have shown the condition of limit cycle and isochronicity for a modified Brusselator model.

The original
four variable reversible $Brusselator$ model\cite{grayscott} which after appropriate elimination of variable results in a simple kinetics of relevant two variables\cite{lav2} as

$$\frac{dx}{dt}=a_1+x^2 y-(\alpha+b) x$$
\begin{equation}
\frac{dy}{dt}=b x-x^2 y
\label{eq31}
\end{equation}
where $x$ and $y$ are the dimensionless concentration of some species. The parameters,  
$a_1,b,\alpha > 0$ follow the properties,  $a_1=\mu(1-\beta) a +\mu \beta$  with $a_1>0\Rightarrow$ either $\mu>0 \hspace{.3 cm}\&\hspace{.3 cm} a<\frac{\beta}{1-\beta}$ or $\mu<0 \hspace{.3 cm}\&\hspace{.3 cm} a>\frac{\beta}{1-\beta}$. Depending on the values of $\beta$  one can choose the conditions accordingly.

Supposing  $z=x+y$ and $u=a_1-\alpha x$ one can transform (\ref{eq31}) into one equation i.e. 
\begin{equation}
\dot{z}=u
\label{eq32}
\end{equation}
Here $x$ and $y$ can be expressed as 
\begin{equation}
x=\frac{a_1-u}{\alpha} \hspace{1cm} \& \hspace{1cm} y=z+\frac{u-a_1}{\alpha}.
\end{equation}

Differentiating (\ref{eq32}) w.r.t. t we get,
$$\ddot{z}=\dot{u}=-\alpha a_1-\alpha x^2 y+\alpha(\alpha+b) x$$
$\Rightarrow$
\begin{equation}
\ddot{z}=(b a_1+\frac{a_1^3}{\alpha^2})-\frac{z a_1^2}{\alpha}+\frac{2 a_1 z \dot{z} }{\alpha}-\dot{z} (\alpha+b+\frac{3 a_1^2}{\alpha^2})+\frac{3 a_1 \dot{z}^2 }{\alpha^2}-\frac{\dot{z}^3}{\alpha^2}-\frac{z \dot{z}^2}{\alpha}.
\label{eq34}
\end{equation}

So if $z_s$ be the stationary point then $\dot{z}$
 and $\ddot{z}$ all are zero which shows from (\ref{eq34}) as,
\begin{equation}
\hspace{3 cm}z_s=\frac{\alpha}{a_1^2}(b a_1+\frac{a_1^3}{\alpha^2}).
\end{equation}

Taking perturbation around the fixed point $z_s$ i.e.  $z=z_s+\xi$ with $\dot{z}=\dot{\xi}$ and $\ddot{z}=\ddot{\xi}$ and substituting this in (\ref{eq34}) and on simplification gives 
\begin{equation}
\ddot{\xi}+F(\xi,\dot{\xi}) \dot{\xi} +G(\xi)=0
\label{eq36}
\end{equation}
where the functions, $F(\xi,\dot{\xi})$ and $G(\xi)$ are $$
F(\xi,\dot{\xi})=-\frac{2 a_1 \xi}{\alpha}-b+\frac{a_1^2}{\alpha^2}+\alpha-\frac{2 a_1 \dot{\xi}}{\alpha^2}+\frac{b \dot{\xi}}{a_1}+\frac{\dot{\xi}^2}{\alpha^2}+\frac{\xi \dot{\xi}}{\alpha};$$
\begin{equation}
\hspace{3 cm}G(\xi)= \frac{a_1^2 \xi}{\alpha}.
\end{equation}
This is $Lienard$ system because all the conditions for being a $Lineard$ oscillator which are stated  previously are satisfied by using the given conditions $a_1,b,\alpha > 0$. 
Thus if there is any stable limit cycle then it must satisfy $F(0,0)<0$ i.e.
\begin{equation}
\hspace{3 cm}a_1^2<(b-\alpha) \alpha^2.
\end{equation}

Now (\ref{eq34}) can be written as,
\begin{equation}
\ddot{\xi}+\omega^2 \xi=\frac{2 a_1 \xi \dot{\xi}}{\alpha}+(b-\frac{a_1^2}{\alpha^2}-\alpha)\dot{\xi}+(\frac{2 a_1}{\alpha^2}-\frac{b}{a_1})\dot{\xi}^2-\frac{\dot{\xi}^3}{\alpha^2}-\frac{\xi \dot{\xi}^2}{\alpha}
\end{equation} 
$$\omega^2=\frac{a_1^2}{\alpha}>0.$$
For book keeping purpose introducing $\lambda(0<\lambda<<1)$ in such a way that the above equation can be expressed as,
\begin{equation}
\ddot{\xi}+\omega^2 \xi=\lambda\frac{2 a_1 \xi \dot{\xi}}{\alpha}+(b-\frac{a_1^2}{\alpha^2}-\alpha)\dot{\xi}+\lambda(\frac{2 a_1}{\alpha^2}-\frac{b}{a_1})\dot{\xi}^2+O(\lambda^2)
\label{eq40}
\end{equation}
(neglecting $O(\lambda^2)$ included terms)

To solve the equation by  applying $RG$ technique, we can take $\xi=\xi_0+\lambda\xi_1+\lambda^2\xi_2+\cdots$ i.e. $\xi=\xi_0+\lambda\xi_1+O(\lambda^2)$, by neglecting $O(\lambda^2)$ terms.
So putting $\xi$ in (\ref{eq40}) and equating the coefficients of both sides for $\lambda^0$ and $\lambda^1$ we  get, 

$\hspace{1.6 cm}\lambda^0:\ddot{\xi_0}+\omega^2 \xi_0=(b-\frac{a_1^2}{\alpha^2}-\alpha)\dot{\xi_0}$
\begin{equation}
\lambda^1:\ddot{\xi_1}+\omega^2 \xi_1=\frac{2 a_1 \xi_0 \dot{\xi_0}}{\alpha}+(b-\frac{a_1^2}{\alpha^2}-\alpha)\dot{\xi_1}+(\frac{2 a_1}{\alpha^2}-\frac{b}{a_1})\dot{\xi_0}^2.
\label{eq41}
\end{equation}

Now if we treat the constant coefficient in $F(\xi,\dot{\xi})$ as zero i.e. $F(0,0)=0$ or, $(b-\frac{a_1^2}{\alpha^2}-\alpha)=0$ then from (\ref{eq41})  we have

$\hspace{1.6 cm}\lambda^0:\ddot{\xi_0}+\omega^2 \xi_0=0$,
\begin{equation}
\lambda^1:\ddot{\xi_1}+\omega^2 \xi_1=\frac{2 a_1 \xi_0 \dot{\xi_0}}{\alpha}+(\frac{2 a_1}{\alpha^2}-\frac{b}{a_1})\dot{\xi_0}^2.
\label{eq42}
\end{equation}

Setting initial condition $\xi(t)=A$ and $\dot{\xi(t)}=0$ at $t=t_0$, then by comparing $\lambda$ we  get $\xi_0=A$ and $\xi_i=0 , \forall i>0$ and $\dot{\xi_i}=0 , \forall i\geq 0$ at $t=t_0$.

After solving (\ref{eq42}), $\xi_0(t)$ and $\xi_1(t)$ becomes,
$$\xi_0(t)=A \cos \omega (t-t_0)$$,

$$\xi_1(t)=\frac{A^2}{2} (\frac{2 a_1}{\alpha^2}-\frac{b}{a_1}) \lbrace1+\frac{\cos 2 \omega (t-t_0)}{3}-\frac{4 \cos \omega (t-t_0)}{3}\rbrace$$
\begin{equation}
\hspace{1.5 cm}+\frac{A^2 a_1}{3 \omega \alpha} \lbrace \sin 2 \omega (t-t_0)-2 \sin \omega (t-t_0) \rbrace.
\end{equation}

So, $\xi(t)$ becomes ( on using $\theta_0=-\omega t_0$),
$$\xi(t)=A \cos (\omega t+\theta_0)+\lambda [\frac{A^2}{2} (\frac{2 a_1}{\alpha^2}-\frac{b}{a_1}) \lbrace1+\frac{\cos 2 (\omega t+\theta_0)}{3}$$
\begin{equation}
\hspace{1.3 cm}-\frac{4 \cos (\omega t+\theta_0)}{3}\rbrace+\frac{A^2 a_1}{3 \omega \alpha} \lbrace \sin 2 (\omega t+\theta_0)-2 \sin (\omega t+\theta_0) \rbrace].
\label{eq44}
\end{equation}
Now considering perturbation in the time interval $(t-t_0)$ by splitting $(t-\tau)+(\tau-t_0)$, where $t_0<\tau<t$ and $\tau$ is very close to $t_0$, we define the interval $(t-\tau)$ as principal part and the remaining small part, $(\tau-t_0)$ can be neglected. 

Suppose here considering perturbation of the time interval the amplitude and the phase slightly change from $A$ to $A(\tau)$ and $\theta_0$ to $\theta(\tau)$. From $RG$ technique the relation between them are $A(\tau)=\frac{A}{Z_1(\tau,t_0)}$  and $\theta(\tau)=\theta_0-Z_2(\tau,t_0)$ where 
\begin{equation}
Z_1(\tau,t_0)=1+\sum_{1}^{\infty} \lambda^n p_n \hspace{1 cm} and \hspace{1 cm} Z_2(\tau,t_0)=0+\sum_{1}^{\infty} \lambda^n q_n.
\end{equation}

Since we are neglecting $O(\lambda^2)$ then from previous equation we get,
\begin{equation}
Z_1(\tau,t_0)=1+\lambda p_1+O(\lambda^2) \hspace{.5 cm} and \hspace{.8 cm} Z_2(\tau,t_0)=\lambda q_1+O(\lambda^2).
\end{equation} 

If we put the functions $Z_1$ and $Z_2$ as well as $A$ and $\theta_0$ in (\ref{eq44}) and remove the terms which could  lead to divergence we must get either $p_1$ is zero or anything containing $(\tau-t_0)$ and finally for $q_1$. But because of the smallness of $(\tau-t_0)$ we can take $p_1$ and $q_1$ approximately zero.
So, after using the above argument the constant, $A$ becomes $A(\tau)$ and constant $\theta_0$ becomes $\theta(\tau)$ i.e. they become dependent upon the time variable $\tau$. Also, if any term   multiplied directly by $(t-t_0)$ in the final solution of $\xi(t)$, then we can convert it into $(t-\tau)$ by neglecting the non-principal part. 


Using all above considerations  in equation (44), $\xi(t)$ becomes,

$$\xi(t)=A(\tau) \cos (\omega t+\theta(\tau))+\lambda [\frac{A^2 (\tau)}{2} (\frac{2 a_1}{\alpha^2}-\frac{b}{a_1}) \lbrace1+\frac{\cos 2 (\omega t+\theta(\tau))}{3}$$
\begin{equation}
\hspace{.8 cm}-\frac{4 \cos (\omega t+\theta(\tau))}{3}\rbrace+\frac{A^2 (\tau) a_1}{3 \omega \alpha} \lbrace \sin 2 (\omega t+\theta(\tau))-2 \sin (\omega t+\theta(\tau)) \rbrace].
\label{eq47}
\end{equation}

So, finally under the condition $(\frac{\partial \xi}{\partial \tau})\vert_{t}=0$, since the final solution can not be dependent on $\tau$, (\ref{eq47}) shows,

$[\cos (\omega t+\theta(\tau))+\lambda A(\tau)(\frac{2 a_1}{\alpha^2}-\frac{b}{a_1}) \lbrace1+\frac{\cos 2 (\omega t+\theta(\tau))}{3}-\frac{4 \cos (\omega t+\theta(\tau))}{3}\rbrace+\frac{2 \lambda A(\tau) a_1}{3 \omega \alpha}\lbrace \sin 2 (\omega t+\theta(\tau))-2 \sin (\omega t+\theta(\tau)) \rbrace]\frac{dA}{d \tau} + [-A(\tau) \sin (\omega t+\theta(\tau))+\frac{\lambda A^2(\tau)}{2 }(\frac{2 a_1}{\alpha^2}-\frac{b}{a_1})\lbrace -\frac{2 \sin 2((\omega t+\theta(\tau))}{3}+\frac{4 \sin (\omega t+\theta(\tau))}{3}\rbrace+\frac{2 \lambda A^2(\tau) a_1}{3 \omega \alpha } \lbrace \cos 2 (\omega t+\theta(\tau))-\cos(\omega t+\theta(\tau))]\frac{d \theta}{d \tau}=0$.

Since, $A (\tau) \neq 0 , a_1 \neq 0 , \alpha \neq 0 , b \neq 0 , \omega \neq0 $, then none of the above which are in third brackets are zero. Therefore the only possible way to balance the equation is $\frac{dA}{d \tau}=0$ and $\frac{d \theta}{d \tau}=0$. 
So, this leads to isochronous oscillator of centre type and finally it cannot have any limit cycle. 




Since, we know that $a_1,b,\alpha>0$ and $a_1=\mu(1-\beta) a +\mu \beta$ so this will give either $\mu>0 \hspace{.3 cm}\&\hspace{.3 cm} a<\frac{\beta}{1-\beta}$ or $\mu<0 \hspace{.3 cm}\&\hspace{.3 cm} a>\frac{\beta}{1-\beta}$. The parameters $a$ and $\beta$ are dependent on each other. Now we deal with the positive region of parameters and  we suppose $\mu=1$ and $a=1$ which gives
 $\beta>0.5$. For this set of parametric values the boundary condition satisfies, $a_1^2=(b-\alpha) \alpha^2$ which produces \ref{fig1a}.
In Figure(\ref{fig1a}) parametric variation of alpha$(\alpha)$ and $b$ are shown in which the boundary line separates the region into stable limit cycle and stable focus centred for the modified Brusselator model.
$Figure$ \ref{fig1} shows a stable limit cycle solution for suitable choice of parameters, $\mu=1$, $a=1$, $\beta=0.6$, $\alpha=2$ and $b=2.5$ together which satisfies the limit cycle condition, $F(0,0)<0$. $Figure$ \ref{fig2} shows a centre type solution satisfying $F(0,0)=0$ by taking suitable choice of parameters, $\mu=1$, $a=1$, $\beta=0.6$, $\alpha=2$ and $b=2.25$. Since $Figure$ \ref{fig2} is a centre type, it must be  closer to the fixed point but not form any limit cycle.

\subsection{Simple Glycolytic Oscillator}

The glycolytic oscillator\cite{gly1,gly2,gly3,glyscott} is mainly observed  in the yeast, which is described with respect to its overall dynamics and biochemical properties  of its enzyme phospho fructokinase. Kinetic properties are complemented by the  mathematical analysis of
  Selkov\cite{gly1} and related models\cite{glyscott}. Here we have considered its modified form for the oscillatory glycolysis in closed vessels by Merkin-Needham-Scott(MNS)\cite{glyscott}  as  
 
$$\dot{x}=-x+(a+x^2) y$$
\begin{equation}
\dot{y}=b-(a+x^2) y.
\label{eq48}
\end{equation}
with $x$ and $y$ corresponding to the intermediate species concentrations.
The phosphofructokinase step considered by Selkov's model
and its MNS-generalization considers ATP to ADP transition accompanying fructose-6-phosphate(F6P) to fructose-1,6-diphosphate(F1,6DP).  
The parameter $b$ means ATP influx and $a$ is the
rate of non-catalyzed side-steps (a side-process, which needs to be taken into
account for the closed vessel consideration, as shown by Merkin-Needham-Scott\cite{glyscott}).
The fixed point of the system is at $x=a,y=\frac{b}{a+b^2}$. It is stable focus for a certain parameter range and an unstable focus for certain others. The crossover from stable to unstable focus occurs on the boundary curve which is a locus of points in the $a-b$ plane where a Hopf bifurcation occurs i.e. the fixed point for those values of $(a,b)$ is a centre which satisfies the equation $(a+b^2)^2 +(a-b^2)=0$ and  can be obtained by checking condition for stability or from the eigenvalues.

If we suppose $z=x+y$ and $u=b-x$ then we can transform (\ref{eq48}) into a  form i.e.
 \begin{equation}
\hspace{4cm}\dot{z}=u.
\label{eq49}
\end{equation}
Here $x , y$ can be expressed as 
\begin{equation}
x=(b-u) \hspace{1cm} \& \hspace{1cm} y=(z-b+u).
\end{equation}
Differentiating (\ref{eq49}) w.r.t. t and eliminating $x$ and $y$ gives,
\begin{equation}
\ddot{z}=-(1+a+3 b^2)\dot{z}-(a+b^2)z+(b+a b+b^3)+2 b z \dot{z}+3 b \dot{z}^2-z \dot{z}^2-\dot{z^3}
\label{eq51}
\end{equation}
If $z_s$ be the fixed point of $z$, for which $\dot{z},\ddot{z}$ all are zero then, \begin{equation}
z_s=b+\frac{b}{a+b^2}=k(say).
\end{equation}

Similarly as in previous case we take a perturbation around $k$, i.e. $z=k+\xi,\dot{z}=\dot{\xi},\ddot{z}=\ddot{\xi}$, then (\ref{eq51}) gives a $Lienard$ system, 
\begin{equation}
\ddot{\xi}+F(\xi,\dot{\xi}) \dot{\xi}+G(\xi)=0
\label{eq53}
\end{equation}
where
$$F(\xi,\dot{\xi})=(1+a+3 b^2)-2 b \xi-2bk-3 b \dot{\xi}+\xi\dot{\xi}+k\dot{\xi}+\dot{\xi}^2$$ and
\begin{equation}
\hspace{3.5 cm} G(\xi)=(a+b^2)\xi.
\end{equation}

Note that all the conditions for a $Lienard$ system are satisfied for suitable choice of $a$ and $b$ which is also obvious for $a\geq0$ whatever $b$ may be. So the condition for existence of stable limit cycle is $F(0,0)<0$ i.e.
\begin{equation}
(a+b^2)+\frac{(a-b^2)}{(a+b^2)}<0.
\end{equation}

Suppose constant coefficient present in $F(\xi,\dot{\xi})$ is taken as zero then (\ref{eq53}) shows, 
\begin{equation}
\ddot{\xi}-\lbrace2 b \xi+3 b \dot{\xi}-\xi\dot{\xi}-k\dot{\xi}-\dot{\xi}^2\rbrace\dot{\xi}
+(a+b^2)\xi=0
\end{equation}
\begin{equation}
\Rightarrow \ddot{\xi}+(a+b^2)\xi=2 b \xi \dot{\xi}+3 b \dot{\xi}^2-\xi\dot{\xi}^2-k\dot{\xi}^2-\dot{\xi}^3,
\end{equation}

Similarly we consider
\begin{equation}
\ddot{\xi}+(a+b^2)\xi=2 \lambda b \xi \dot{\xi}+3\lambda b \dot{\xi}^2-\lambda^2\xi\dot{\xi}^2-k\lambda\dot{\xi}^2-\lambda^2\dot{\xi}^3
\end{equation}
where $\lambda(0<\lambda<<1)$. So, because of smallness of $\lambda$, neglecting $O(\lambda^2)$,
one finds
\begin{equation}
\ddot{\xi}+(a+b^2)\xi=2 \lambda b \xi \dot{\xi}+3\lambda b \dot{\xi}^2-k\lambda\dot{\xi}^2+O(\lambda^2).
\label{eq59}
\end{equation}

Taking a perturbative solution of $\xi$ as $\xi=\xi_0+\lambda\xi_1+\lambda^2\xi_2+\lambda^3\xi_3+\cdots
$ i.e. $\xi=\xi_0+\lambda\xi_1+O(\lambda^2)$ (neglecting $O(\lambda^2)$) and using above perturbative solution and comparing the coefficients of $\lambda^0,\lambda^1$ (\ref{eq59}) gives,

$\hspace{1.7 cm}\lambda^0:\ddot{\xi_0}+(a+b^2)\xi_0=0$
\begin{equation}
\lambda^1:\ddot{\xi_1}+(a+b^2)\xi_1=2 b \xi_0 \dot{\xi_0}+(3 b - k)\dot{\xi_0}^2.
\label{eq60}
\end{equation}
Using the most general initial condition $\xi(t)=A$ and $\dot{\xi(t)}=0$ at $t=t_0$, then by comparing $\lambda$ as similar as in previous case we must get $\xi_0=A$ and $\xi_i=0 , \forall i>0$ with $\dot{\xi_i}=0 , \forall i\geq0$ at $t=t_0$.

After solving (\ref{eq60}), $\xi_0(t)$ and $\xi_1(t)$ becomes,
$$\xi_0(t)=A \cos \omega (t-t_0)$$ and 
$\xi_1(t)=-\frac{2bA^2}{3 \omega}\sin \omega(t-t_0)-\frac{2(3b-k)A^2}{3}\cos \omega(t-t_0)+\frac{(3b-k)A^2}{2}\lbrace 1+$
\begin{equation}
\frac{\cos 2\omega(t-t_0)}{3} \rbrace+\frac{bA^2}{3 \omega}\sin 2\omega(t-t_0)
\end{equation}
where $\hspace{.5 cm}\omega^2=(a+b^2)>0$.

So $\xi(t)$ becomes,

$\xi(t)=A \cos \omega (t-t_0)+\lambda[-\frac{2bA^2}{3 \omega}\sin \omega(t-t_0)-\frac{2(3b-k)A^2}{3}\cos \omega(t-t_0)$
\begin{equation}
+\frac{(3b-k)A^2}{2}\lbrace 1+
\frac{\cos 2\omega(t-t_0)}{3} \rbrace+\frac{bA^2}{3 \omega}\sin 2\omega(t-t_0)].
\end{equation}
If $\theta_0=-\omega t_0$, then the above equation  can be written as,

$\xi(t)=A \cos (\omega t+\theta_0)+\lambda[-\frac{2bA^2}{3 \omega}\sin (\omega t+\theta_0)-\frac{2(3b-k)A^2}{3}\cos (\omega t+\theta_0)$
\begin{equation}
+\frac{(3b-k)A^2}{2}\lbrace 1+
\frac{\cos 2(\omega t+\theta_0)}{3} \rbrace+\frac{bA^2}{3 \omega}\sin 2(\omega t+\theta_0)].
\label{eq63}
\end{equation}

Now adding another perturbation in the time interval $(t-t_0)$ by splitting $(t-\tau)+(\tau-t_0)$, where $t_0<\tau<t$ and $\tau$ is very close to $t_0$ we define the interval $(t-\tau)$ as a principal part and the remaining part $(\tau-t_0)$ can be neglected because of the smallness. 
Suppose that on taking perturbation the time interval the amplitude and the phase slightly be changed from $A$ to $A(\tau)$ and $\theta_0$ to $\theta(\tau)$. From $RG$ technique the relation between them are $A(\tau)=\frac{A}{Z_1(\tau,t_0)}$, and $\theta(\tau)=\theta_0-Z_2(\tau,t_0)$; where 
\begin{equation}
Z_1(\tau,t_0)=1+\sum_{1}^{\infty} \lambda^n p_n \hspace{1 cm} and \hspace{1 cm} Z_2(\tau,t_0)=0+\sum_{1}^{\infty} \lambda^n q_n.
\end{equation}
Since we are neglecting $O(\lambda^2)$ then from previous equation we get,
\begin{equation}
Z_1(\tau,t_0)=1+\lambda p_1+O(\lambda^2) \hspace{.5 cm} and \hspace{.8 cm} Z_2(\tau,t_0)=\lambda q_1+O(\lambda^2).
\end{equation} 

Now, if we put the functions $Z_1$ and $Z_2$ as well as $A$ and $\theta_0$ in (\ref{eq63}) and remove the terms which could  lead to divergence we must get either $p_1$ is zero or anything containing $(\tau-t_0)$ and same for $q_1$. But because of the smallness of $(\tau-t_0)$ we can take $p_1$ and $q_1$ approximately to zero.


Using all above results in the equation (\ref{eq63}) we get,

$\xi(t)=A(\tau) \cos (\omega t+\theta(\tau))+\lambda[-\frac{2bA^2(\tau)}{3 \omega}\sin (\omega t+\theta(\tau))-\frac{2(3b-k)A^2(\tau)}{3}\cos (\omega t+\theta(\tau))$
\begin{equation}
+\frac{(3b-k)A^2(\tau)}{2}\lbrace 1+
\frac{\cos 2(\omega t+\theta(\tau))}{3} \rbrace+\frac{bA^2(\tau)}{3 \omega}\sin 2(\omega t+\theta(\tau))]
\label{eq66}
\end{equation}

So finally under the condition as in $RG$ method $(\frac{\partial \xi}{\partial \tau})\vert_{t}=0$ (\ref{eq66}) gives ,
$$\frac{dA}{d \tau}=0 \hspace{1 cm} and$$ $$\hspace{-1.35 cm}\frac{d \theta}{d \tau}=0.$$
Since, $a \neq 0 , b \neq 0 (\Rightarrow  \omega \neq 0)$ and $A (\tau) \neq 0$ then none of the above in brackets in the last equation are zero which are obtained from $RG$ condition. Therefore the only possible way is  balancing the equation by making them zero. So this leads to isochronous centre type and finally it cannot have any limit cycle. So we can construct an isochronous oscillatory equation from $Lienard$ system by suitable choice of parameters which makes $F(0,0)=0$.

$Figure$ \ref{fig4} gives a stable limit cycle for describing glycolytic oscillator model when $a$ and $b$ are chosen as 0.11 and 0.6, respectively, together satisfies limit cycle condition. $Figure$ \ref{fig5} represents a centre type solution when $a$ is fixed at zero and and $b$ is fixed at 1 (which are  on the boundary point of $Figure$ \ref{fig3} together which satisfies $F(0,0)=0$).

\subsection{van der Pol Type Oscillator Model}

Here we take  an example of  $van$ $der$ $Pol$ type oscillator. Van der Pol oscillator\cite{murraymbio,strogatz,chen1,len3.5} arises in many nonlinear dynamical systems including in chemical oscillation with cubic nonlinear processes. This case is readily convertable to a Lienard oscillator form. Here we show the conditions of limit cycle and isochronicity with a slightly different analysis than the previous examples. The set of differential equation for $van$ $der$ $Pol$ type oscillator is given as, 

$$\dot{x}=y$$ 
\begin{equation}
\hspace{2.5 cm}\dot{y}=-\epsilon y (x^2-a^2)-\omega^2 x, a\in \Re.
\label{eq67}
\end{equation}
Originally in the $van$ $der$ $Pol$ oscillator equation $a$ is  1. The size of the limit cycle depends on the magnitude of a. Here our analysis is valid for $a\in \Re$, however, for the purpose of RG analysis\cite{chen1,chen2} here we consider a little generalized form with $\epsilon$ as a smallness parameter of perturbation although fixed by the van der Pol system.
There is a fixed point at the origin which has a stable focus for $\epsilon<0$ and unstable focus for $\epsilon>0$. The fixed point is a centre for $\epsilon=0$. If we differentiate first equation w.r.t. $t$ and use second equation we get 
\begin{equation}
\ddot{x}+\epsilon \dot{x} (x^2-a^2)+\omega^2 x=0.
\label{eq68}
\end{equation}
It is of $Lienard$ form and gives a $van$ $der$ $Pol$ oscillator for $a=1$. We now  analyse the model without taking any particular value of $a$. For this model the damping force is $F(x,\dot{x})=\epsilon(x^2-a^2)$ and $G(x)=\omega^2 x$. So they satisfy all the conditions for being a $Lienard$ system which can be a limit cycle if $F(0,0)<0$ i.e. $\forall a\in \Re_{\ne 0} \hspace{.3cm} and \hspace{.3cm} \epsilon>0$ ($\epsilon$ is very  small). If  $\epsilon<0$ then $a^2<0$, which can not be possible for any real $a$. Now for further calculation we rewrite the equation (\ref{eq68}) as, 

\begin{equation}
\ddot{x}+\omega^2 x=-\epsilon \dot{x} (x^2-a^2).
\label{eq69}
\end{equation}

Expanding $x(t)$ as $x(t)=x_0(t)+\epsilon x_1(t)+\epsilon^2 x_2(t)+\cdots$, i.e. $x(t)=x_0(t)+\epsilon x_1(t)+O(\epsilon^2)$, then neglecting $O(\epsilon^2)$ $(0<\epsilon<<1)$ and from (\ref{eq69}) by comparing the coefficients of $\epsilon^0,\epsilon^1$ we get,
$$\epsilon^0:\ddot{x_0}+\omega^2 x_0=0$$
\begin{equation}
\epsilon^1:\ddot{x_1}+\omega^2 x_1=-\dot{x_0} (x_0^2-a^2). 
\label{eq70}
\end{equation}

Using the general initial condition, $x(t)=A$ and $\dot{x(t)}=0$ at $t=t_0$, then by comparing $\epsilon$ as similar in previous case we get $x_0=A$ and $x_i=0 , \forall i>0$ along with $\dot{x_i}=0 , \forall i\geq0$ at $t=t_0$. After solving (\ref{eq70}), $x_0(t)$ and $x_1(t)$ becomes,
$$x_0(t)=A \cos \omega (t-t_0)$$ and
$$x_1(t)=\frac{(7 A^3-16 A a^2)}{32 \omega} \sin \omega(t-t_0)-\frac{(A^3-4 A a^2)}{8}(t-t_0)\cos \omega (t-t_0)$$
\begin{equation}
-\frac{A^3}{32 \omega} \sin 3 \omega (t-t_0).
\end{equation}
Then $x(t)$ becomes on using $\theta_0=-\omega t_0$,

$$x(t)=A \cos (\omega t+\theta_0)+\epsilon[\frac{(7 A^3-16 A a^2)}{32 \omega} \sin (\omega t+\theta_0)$$
\begin{equation}
-\frac{(A^3-4 A a^2)}{8}(t-t_0)\cos (\omega t+\theta_0)-\frac{A^3}{32 \omega} \sin 3 (\omega t+\theta_0)].
\label{eq72}
\end{equation}

Now considering  some perturbation in the time interval, $(t-t_0)$ by splitting $(t-\tau)+(\tau-t_0)$, where $t_0<\tau<t$ and $\tau$ is very close to $t_0$. We define the interval $(t-\tau)$ as a principal part and the remaining part $(\tau-t_0)$ can be neglected because of smallness. 

Suppose that on taking perturbation on the time interval the amplitude and the phase slightly be changed from $A$ to $A(\tau)$ and $\theta_0$ to $\theta(\tau)$. Then from $RG$ technique the relation between them are obtained as $A(\tau)=\frac{A}{Z_1(\tau,t_0)}$  and $\theta(\tau)=\theta_0-Z_2(\tau,t_0)$ where 
\begin{equation}
Z_1(\tau,t_0)=1+\sum_{1}^{\infty} \epsilon^n p_n \hspace{1 cm} and \hspace{1 cm} Z_2(\tau,t_0)=0+\sum_{1}^{\infty} \epsilon^n q_n.
\end{equation}
Since we are neglecting $O(\epsilon^2)$ then we obtain,
\begin{equation}
Z_1(\tau,t_0)=1+\epsilon p_1+O(\epsilon^2) \hspace{.5 cm} and \hspace{.8 cm} Z_2(\tau,t_0)=\epsilon q_1+O(\epsilon^2).
\end{equation} 

Now if we put the functions $Z_1$ and $Z_2$ as well as $A$ and $\theta_0$ in (\ref{eq72}) and remove the terms which could lead to divergence we must get either $p_1$ is zero or anything containing $(\tau-t_0)$ and the same for $q_1$. But because of the smallness of $(\tau-t_0)$ we can take $p_1$ and $q_1$ approximately to zero.
So after using the above consideration the constant $A$ becomes $A(\tau)$ and constant $\theta_0$ becomes $\theta(\tau)$ i.e. they become dependent upon the time variable, $\tau$. Again if any term  is multiplied directly by $(t-t_0)$ in the final solution of $x(t)$ then it can be converted to $(t-\tau)$ by neglecting the non-principal part. Here a term is present which is directly multiplied with $(t-t_0)$, so neglecting it and using all above results in the equation (\ref{eq72}) we get,

\hspace*{1 cm}$$x(t)=A(\tau) \cos (\omega t+\theta(\tau))+\epsilon[\frac{\lbrace 7 A^3(\tau)-16 A(\tau) a^2\rbrace}{32 \omega} \sin (\omega t+\theta(\tau))$$
\begin{equation}
-\frac{\lbrace A^3(\tau)-4 A(\tau) a^2\rbrace}{8}(t-\tau)\cos (\omega t+\theta(\tau))-\frac{A^3(\tau)}{32 \omega} \sin 3 (\omega t+\theta(\tau))].
\end{equation}

Finally $\frac{\partial x}{\partial\tau}\mid_t =0$ gives,

$[\cos (\omega t+\theta(\tau))+ \epsilon \lbrace \frac{(21 A^2(\tau)-16 a^2) }{32 \omega}\sin(\omega t+\theta(\tau))-\frac{(3 A^2(\tau)-4 a^2)}{8} (t-\tau)\cos(\omega t+\theta(\tau))-\frac{3 A^2(\tau)}{32 \omega} \sin 3(\omega t+\theta(\tau))\rbrace]\frac{dA}{d\tau}+[-A(\tau)\sin (\omega t+\theta(\tau))+\epsilon \lbrace \frac{(7 A^3(\tau)-16 a^2 A(\tau)) }{32 \omega}\cos(\omega t+\theta(\tau))+\frac{(A^3(\tau)-4 a^2 A(\tau))}{8} (t-\tau)\sin(\omega t+\theta(\tau))-\frac{3 A^3(\tau)}{32 \omega} \cos 3(\omega t+\theta(\tau))\rbrace]\frac{d \theta}{d \tau}=-\epsilon \frac{(A^3(\tau)-4 a^2 A(\tau))}{8} \cos(\omega t+\theta(\tau))$

which leads to $$\hspace{1.4 cm}\frac{dA}{d\tau}=\frac{\epsilon A(\tau)}{2}\lbrace a^2- \frac{ A^2(\tau)}{4}\rbrace+O(\epsilon^2)\hspace{1 cm} and$$ $$\hspace{-1.5cm}\frac{d \theta}{d \tau}=0.$$

Now, if we convert the limit cycle condition, $- \epsilon a^2 < 0$  to $- \epsilon a^2 = 0$ i.e. $a=0$ then it reduces to $\frac{dA}{d\tau}=-\frac{\epsilon A^3(\tau)}{8}+O(\epsilon^2)\hspace{.1 cm} and \hspace{.1cm}\frac{d \theta}{d \tau}=0$, showing an isochronous centre type solution. It also satisfies the solution of $van$ $der$ $Pol$ oscillation, $a=1$ analysed in \cite{len3.5}. For $a=0$, the oscillatory equation also reduces to a centre type as $\epsilon \rightarrow 0$ in the $van$ $der$ $Pol$ type system.

$Figure$ \ref{fig7} shows a stable limit cycle oscillation when $a=0.5$ and $Figure$ \ref{fig8} shows a centre type oscillation when $a$ is fixed at zero. So we obtain a centre type oscillator  which is isochronous. Thus for $F(0,0)=0$, limit cycle  breaks down and it gives an isochronous oscillator. It is found numerically that a limit cycle can be obtained in this system for generalized $van$ $der$ $Pol$ system, $a \in \Re_{\neq 0}$ when  $0<\epsilon a^2<8.14$ satisfies.


\section{Conclusion}

By casting  a class of chemical  oscillations usually governed by two-variable kinetic equations into the form of a $Lienard$ oscillator 
here we have found the conditions of limit cycle and isochronicity.
It is shown that  the conditions are dictated by the 
nonlinear damping coefficient
and the potential or the forcing term which can be controlled by the suitable choice of 
 the experimental parameters of the  chemical oscillators. Although the conditions of limit cycle and isochronicity are shown here with two variables, this mathematical method along with its numerical applicability can also be important for real  higher order system. 
More specifically the main findings in this work are as follows.

$\bullet$  When the two dimensional kinetic equations are transformed into a $Lienard$
system of equation the condition for limit cycle and isochronicity can be stated in a
unified way. In terms of the $Lienard$ type oscillator, the condition for limit cycle is
given by $F (0, 0) < 0$ whereas for the condition of satisfying an isochronous oscillator is
$F (0, 0) = 0$.

$\bullet$ When the limit cycle condition i.e. $F(0,0)<0$ modifies to its boundary i.e. $F(0,0)=0$ depending on the suitable choice of parameters, the $Lienard$ type oscillator transforms into an isochronous oscillator.

$\bullet$ For any $Lienard$ system, when it converts into an isochronous oscillator the system looses its limit cycle stability and it becomes  of centre type.

\noindent
{\bf Acknowledgement}

\noindent
{ Sandip Saha acknowledges RGNF, UGC, India for the partial financial support.}

\newpage

\noindent
{\large{\bf{Figure Captions}}}

\begin{figure}[tbh]
\rotatebox{0}
{\includegraphics[width=7cm,keepaspectratio]{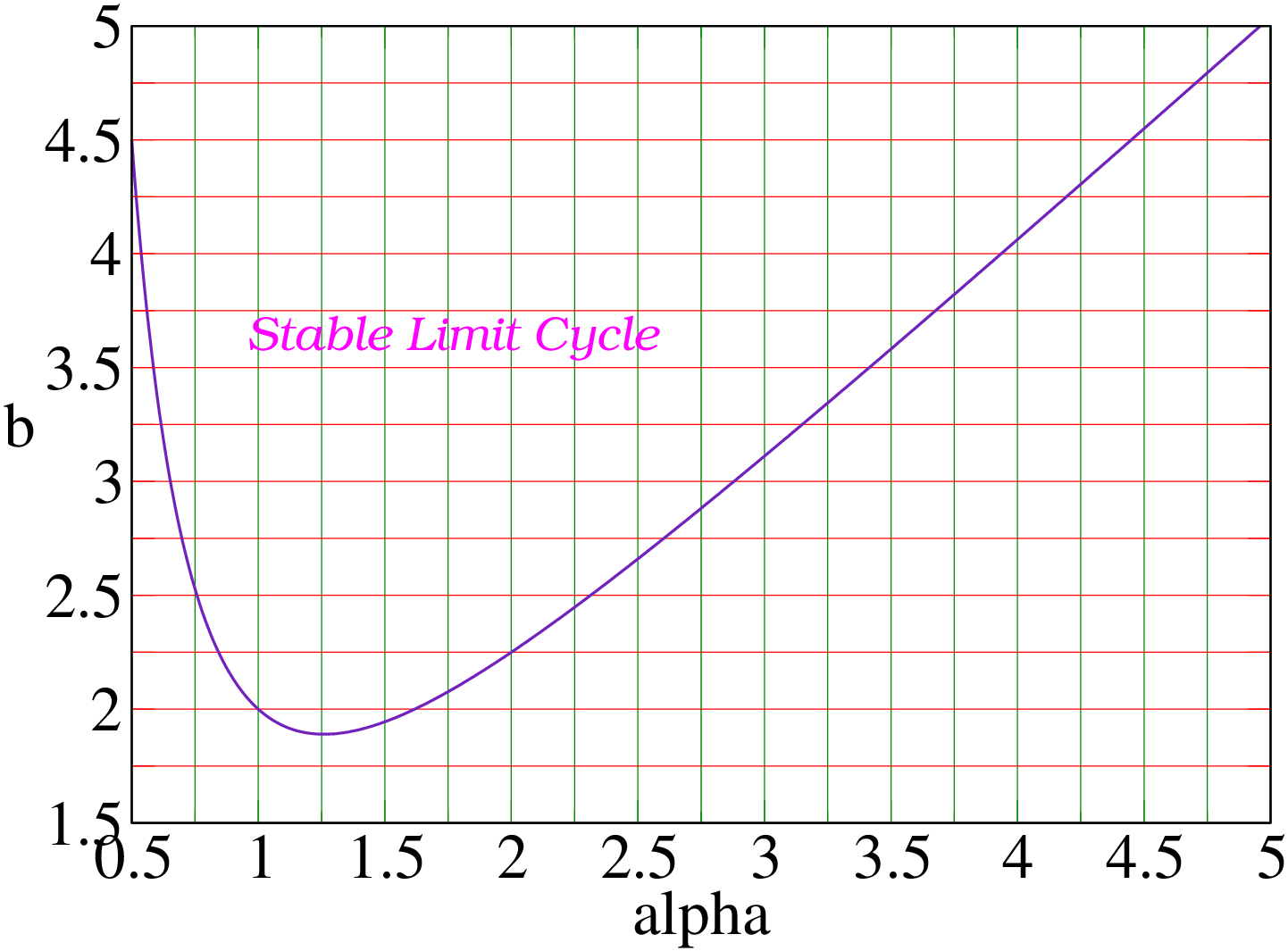}}
\centering
\caption{}
\label{fig1a}
\end{figure}

\ref{fig1a}$Modified$ $Brusselator$ $Model:$ Parametric space diagram for $alpha(\alpha)$ and $b$ in which the boundary line separates the region into stable limit cycle and stable focus when $\mu=1$, $a=1$ and $\beta=0.6$.

\begin{figure}[tbh]
\rotatebox{0}
{\includegraphics[width=7cm,keepaspectratio]{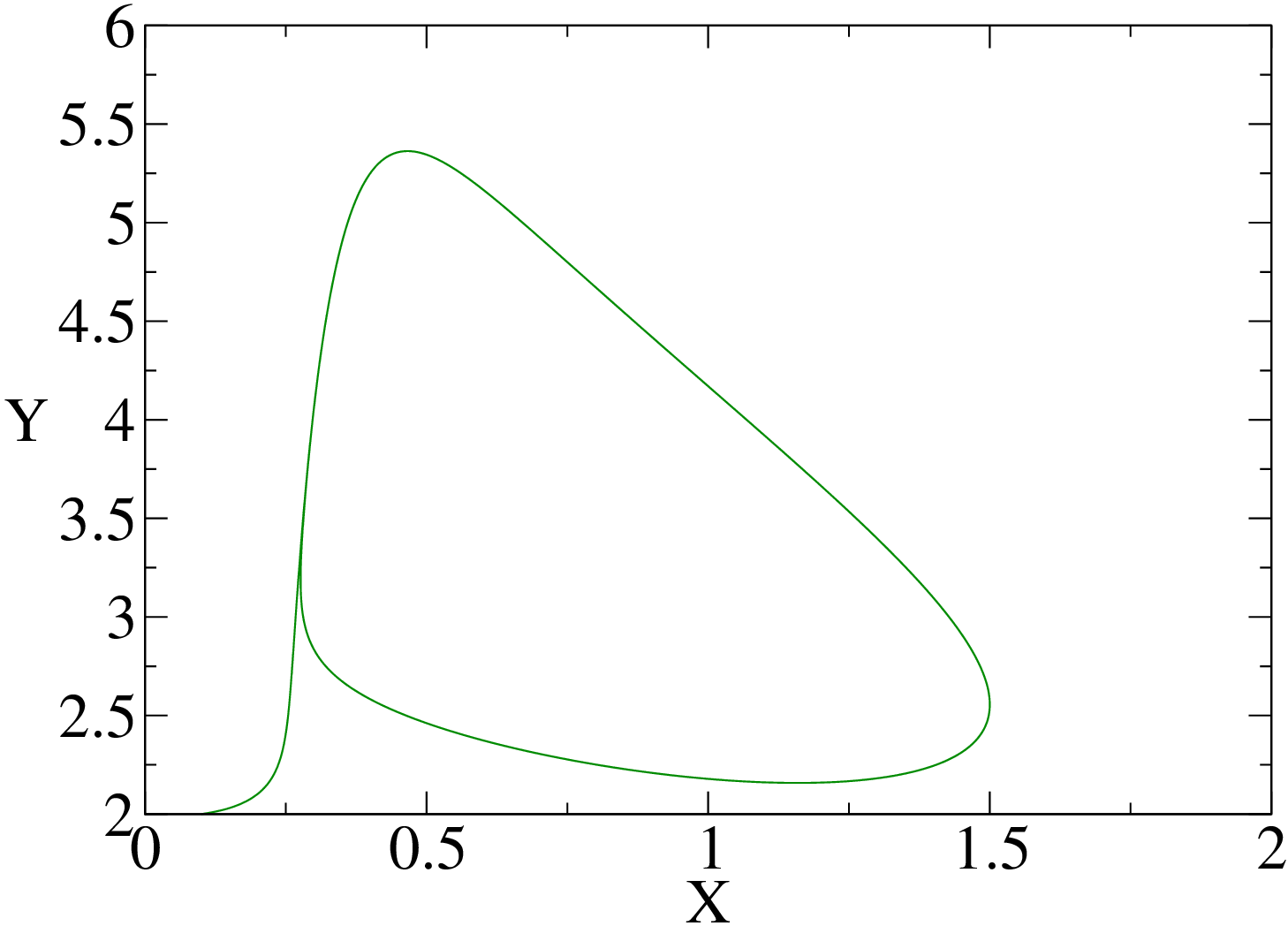}}
\centering
\caption{}
\label{fig1}
\end{figure}

\ref{fig1}$Modified$ $Brusselator$ $Model:$ Phase portrait of (\ref{eq31}) gives a stable limit cycle for suitable choice of parameters, $\mu=1$, $a=1$, $\beta=0.6$, $alpha=2$, $b=2.5$ together which satisfies the limit cycle condition, $F(0,0)<0$.

\begin{figure}[tbh]
\rotatebox{0}
{\includegraphics[width=7cm,keepaspectratio]{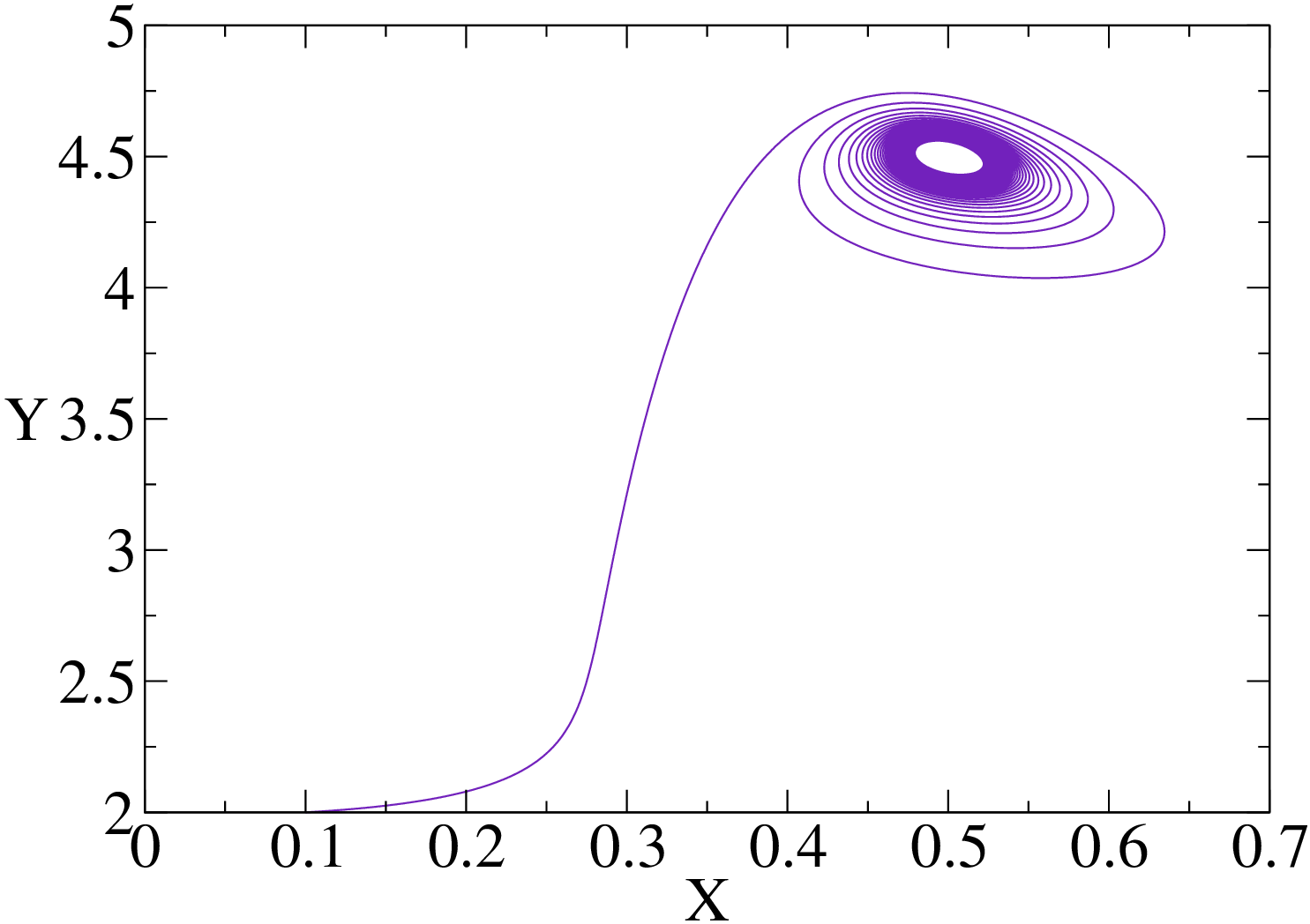}}
\centering
\caption{}
\label{fig2}
\end{figure}

\ref{fig2}$Modified$ $Brusselator$ $Model:$ Phase portrait of (\ref{eq31}) gives a centre type solution when $F(0,0)=0$ by taking suitable choice of parameters, $\mu=1$, $a=1$, $\beta=0.6$, $alpha=2$, $b=2.25$.

\begin{figure}[tbh]
\rotatebox{0}
{\includegraphics[width=7cm,keepaspectratio]{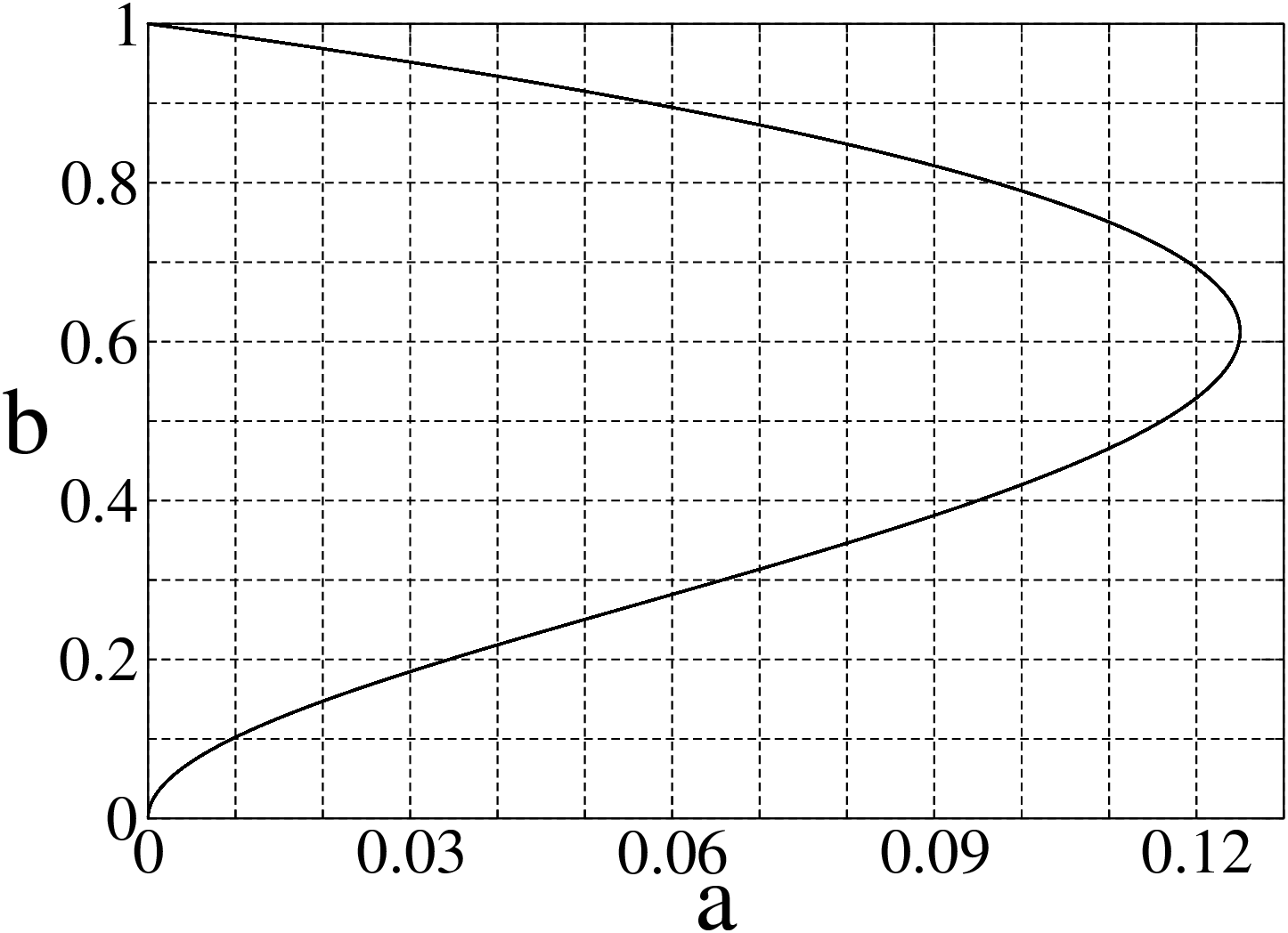}}
\centering
\caption{}
\label{fig3}
\end{figure}

\ref{fig3}$Glycolytic$ $Oscillator:$ Parametric phase portrait for the parameters, $a$ and $b$ in which the boundary line separates the region into stable limit cycle(inner side) and stable focus(outer side).

\begin{figure}[tbh]
\rotatebox{0}
{\includegraphics[width=7cm,keepaspectratio]{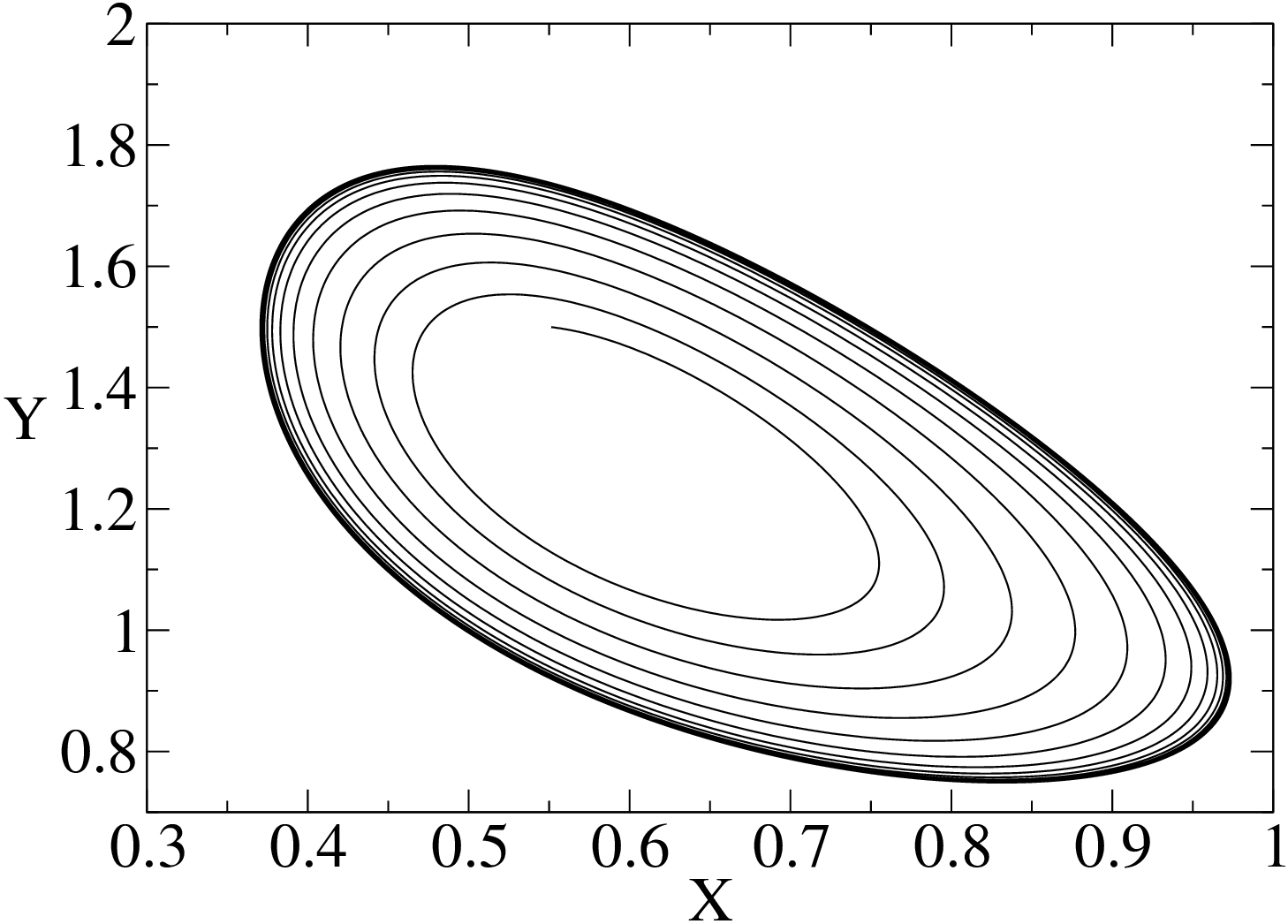}}
\centering
\caption{}
\label{fig4}
\end{figure}

\ref{fig4}$Glycolytic$ $Oscillator:$ Phase space diagram of (\ref{eq48}) when a=0.11 and b=0.6 satisfies limit cycle condition, $F(0,0)<0$ and gives a stable limit cycle.

\begin{figure}[tbh]
\rotatebox{0}
{\includegraphics[width=7cm,keepaspectratio]{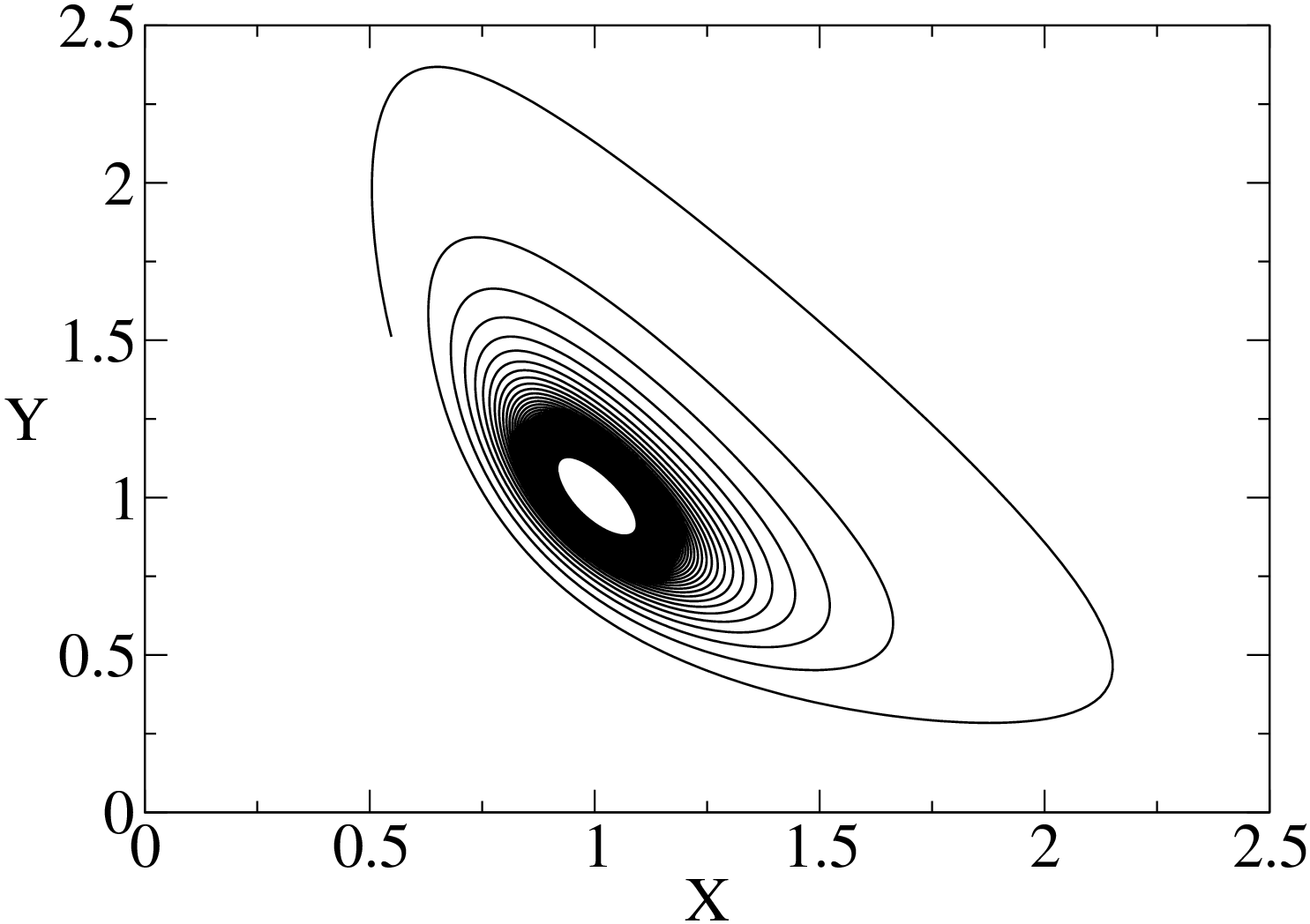}}
\centering
\caption{}
\label{fig5}
\end{figure}

\ref{fig5}$Glycolytic$ $Oscillator:$ Phase space diagram of (\ref{eq48}) when $a=0$ and $b=1$ together which satisfies $F(0,0)=0$ and gives a centre i.e. $a,b$ are  on the boundary line of $Figure$ \ref{fig3}.

\begin{figure}[tbh]
\rotatebox{0}
{\includegraphics[width=7cm,keepaspectratio]{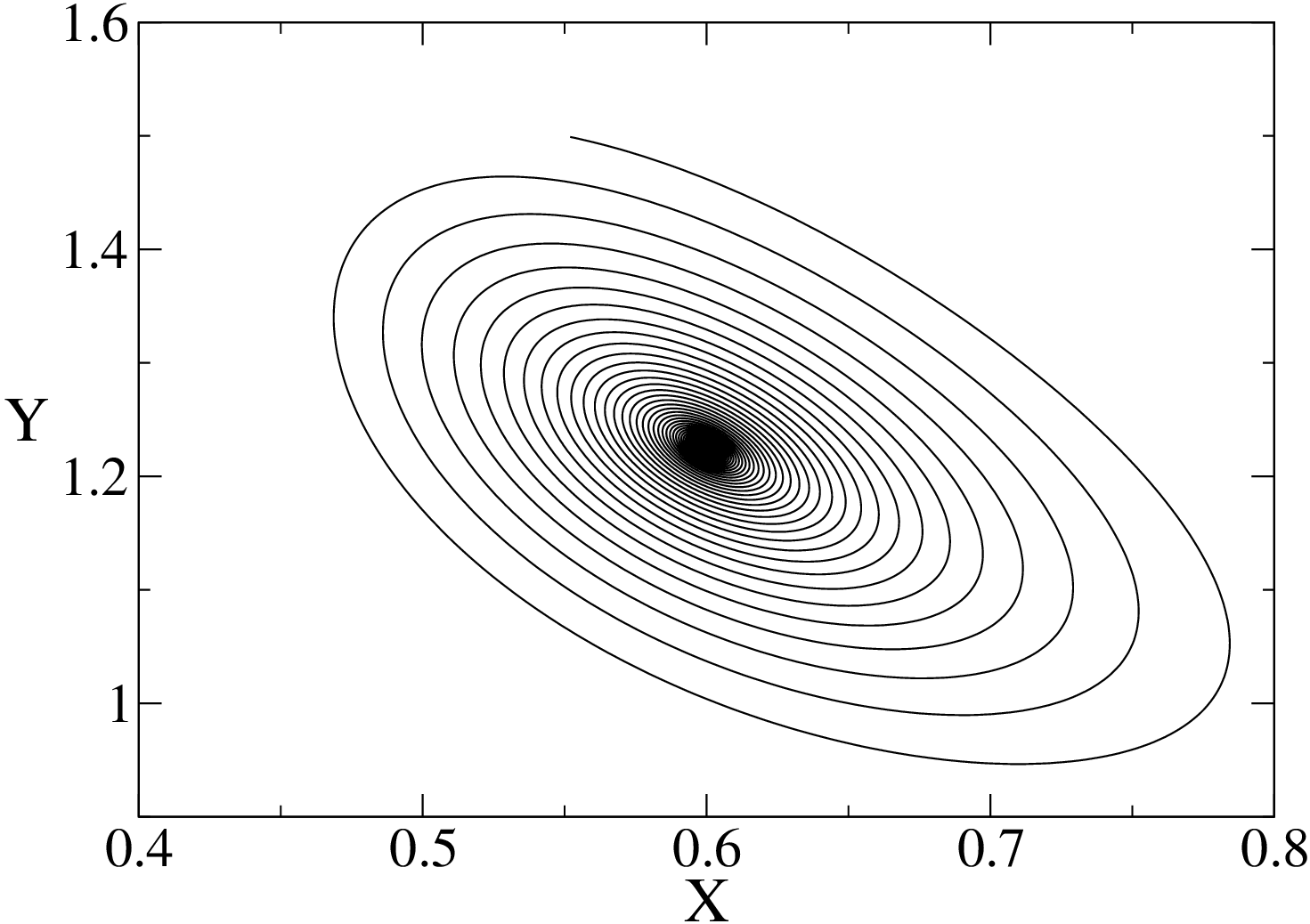}}
\centering
\caption{}
\label{fig6}
\end{figure}

\ref{fig6}
$Glycolytic$ $Oscillator:$ Phase portrait of (\ref{eq48}) when parameters are chosen in such a way that $F(0,0)>0$ i.e. a=0.13 and b=0.6  lie on stable focus.

\begin{figure}[tbh]
\rotatebox{0}
{\includegraphics[width=7cm,keepaspectratio]{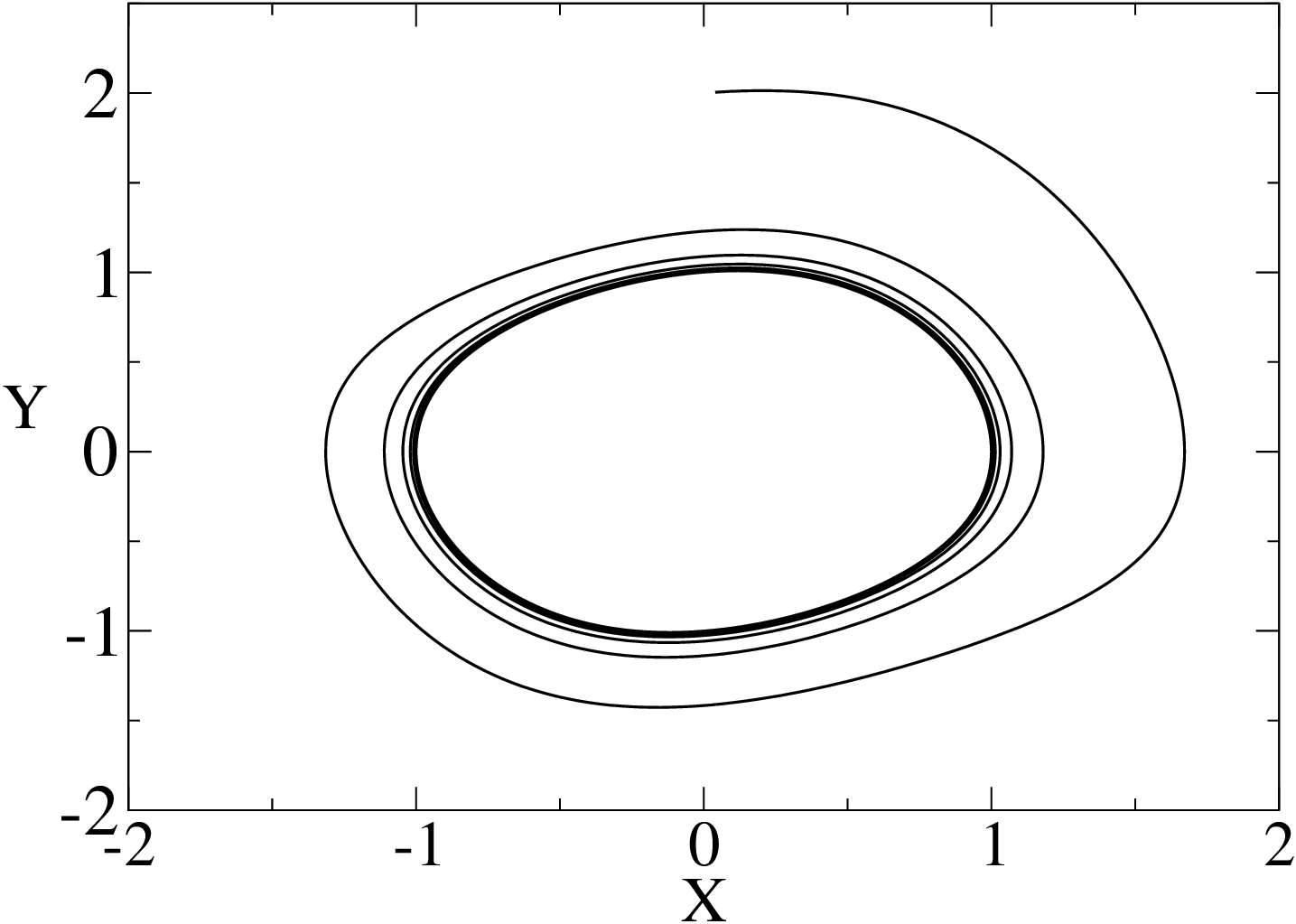}}
\centering
\caption{}
\label{fig7}
\end{figure}

\ref{fig7}$van$ $der$ $Pol$ $Oscillator:$ Phase portrait of (\ref{eq67}) when $a=0.5$ gives a stable limit cycle for $F(0,0)<0$.

\begin{figure}[tbh]
\rotatebox{0}
{\includegraphics[width=7cm,keepaspectratio]{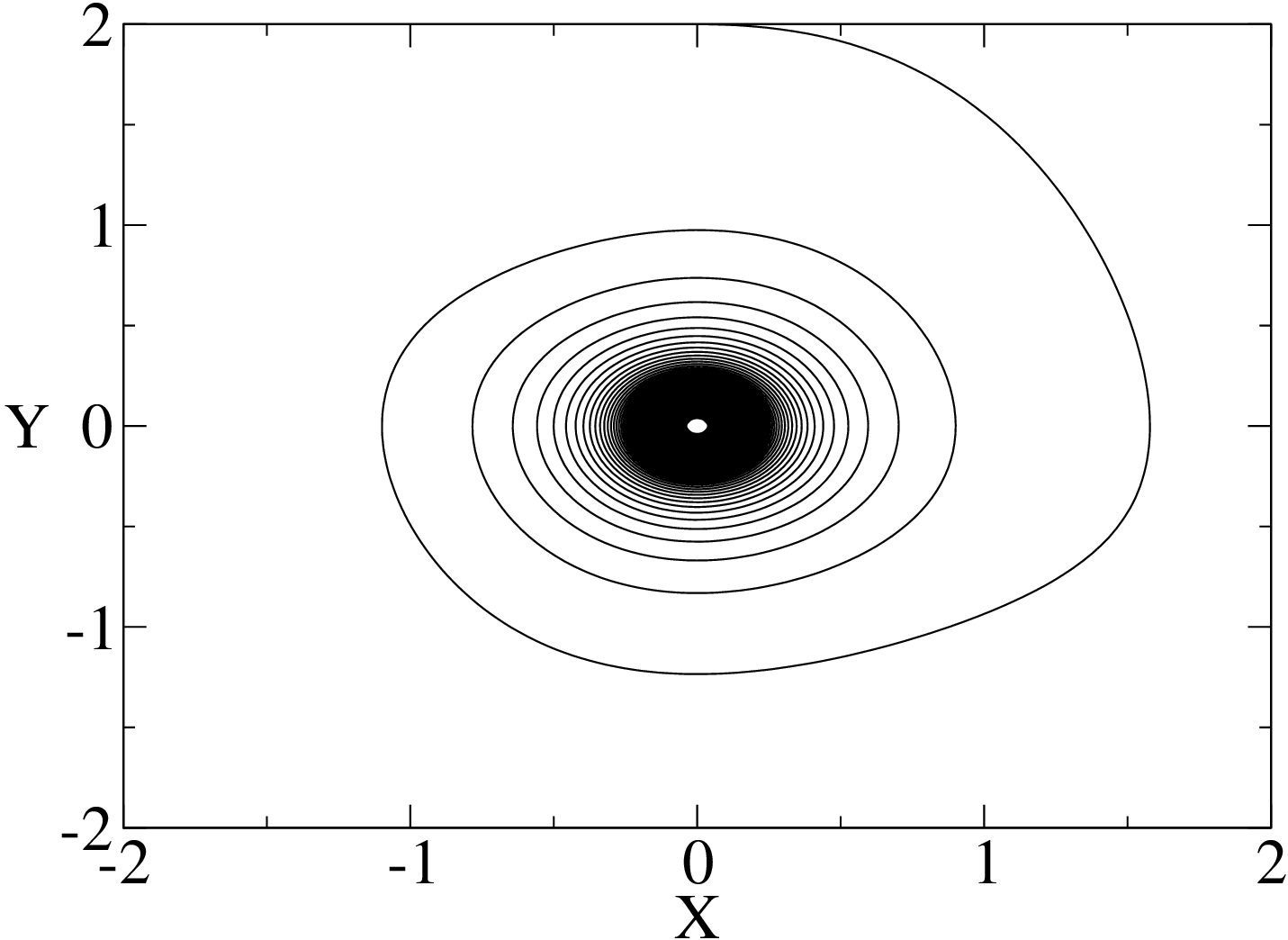}}
\centering
\caption{}
\label{fig8}
\end{figure}

\ref{fig8} $van$ $der$ $Pol$ $Oscillator:$ Phase space diagram of (\ref{eq67}) when $a=0$ satisfies $F(0,0)=0$ which leads to a centre.

\newpage
\newpage

\newpage

\newpage

\newpage

\newpage
\newpage

\newpage

\newpage

\newpage

\end{document}